\documentclass[10pt]{amsart}	
\usepackage{amsmath,amsfonts,amssymb}
\usepackage{verbatim}

\numberwithin{equation}{section}

   \theoremstyle{plain}
\newtheorem{theorem}{Theorem}[section]
\newtheorem{proposition}[theorem]{Proposition}
\newtheorem{lemma}[theorem]{Lemma}
\newtheorem{corollary}[theorem]{Corollary}

	 \theoremstyle{definition}
\newtheorem{definition}[theorem]{Definition}
\newtheorem{example}[theorem]{Example}

	 \theoremstyle{remark}
\newtheorem{remark}{Remark}

\newcommand\tiff{\mbox{ iff }}
\newcommand\tif{\mbox{ if }}
\newcommand\tthen{\mbox{ then }}

\newcommand\tand{\mbox{ and }}

\newcommand\tforall{\mbox{ for all }}

\newcommand\mb[1]{\mbox{ #1 }}

\newcommand\eq\Leftrightarrow
\newcommand\im\Rightarrow

\newcommand\deq{\hfill$\lhd$}

\newcommand\ESP{\textsf{ESP}}
\newcommand\SP{\textsf{SP}}
\newcommand\EWR{\textsf{EWR}}
\newcommand\nrm{\mathrm{nrm}}
\renewcommand\j{\mathrm{j}}

\newenvironment{myenum}{\begin{enumerate}

        \setlength{\itemsep}{0pt}}
        {\end{enumerate}}
\newcommand\thitem[1]{\item[{\upshape (#1)}]}

\title[Sectional pseudocomplementation in posets]{Extension of sectional pseudocomplementation in posets}
\author[J.\ C\={\i}rulis]{J\=anis C{\=\i}rulis}
\address{Institute of Mathematics and Computer Science, University of Latvia,
Raina b., 29,  R\={\i}ga, LV-1050, Latvia}
\email{jc@lanet.lv}

\begin{document}

\begin{abstract}
Sectional pseudocomplementation (sp-comple\-men\-ta\-tion) on a po\-set is a partial operation $*$ which associates with every pair $(x,y)$ of  elements, where $x \ge y$, the pseudocomplement $x*y$ of $x$ in the upper section $[y)$. Any total extension $\to$ of $*$ is said to be an extended sp-complementation and is considered as an implication-like operation. Extended sp-complementations have already be studied on semilattices and lattices. We  describe several naturally arising classes of general posets with extended sp-complementation, present respective elementary properties of this operation,  demonstrate that two other known attempts to isolate particular such classes are in fact not quite correct, and suggest suitable improvements. 
\end{abstract}

\maketitle

\section{Introduction}  \label{intro}

Let $P$ be a poset. If it has the least element $0$, then a pseudocomplement of an element $a$ is, by definition, the greatest element $u$ disjoint with $x$, i.e., such that $0$ is the single lower bound of $u$ and $a$. If every element of $P$ has a pseudocomplement, the poset is said to be pseudocomplemented. Pseudocomplemented posets were introduced in \cite{K2,K3} and independently studied also in \cite{V} and \cite{H}.

Pseudocomplements can be considered also in subsets of $P$ having the least element. Sectionally pseudocomplemneted meet semilattices, i.e., meet semilattices with zero and pseudocomplemented lower sections (principal order ideals) have been considered as early as in \cite{K1,K4}. Semilattices with pseudocomplemented segments also have been called sectionally pseudocomplemented; see, e.g.,  \cite{K5,K6}). In the current literature, this attribute  usually concerns posets with pseudocomplemented upper sections (principal order filters). Such posets with a greatest, or unit element 1  have been discussed by the author in \cite{C-impl,C-annih}; however, join (meet) semilattices with pseudocomplemented principal filters appeared already in \cite{A1} (respectively,  \cite{Sch}), and  lattices, in \cite{Ch-ext}. The term `sectionally pseudocomplemented' is used in this paper in the latter sense. (A variant of this meaning of the term should be mentioned: in some papers, e.g., \cite{ChLP,ChR}, a lattice is said to be sectionally pseudocomplemented if the pseudocomplement of $a \vee b$ in the principal filter $[a)$ exists for all $a,b$. This view on sectional complementedness has been quite recently generalized in \cite{ChL,ChLP} for posets.)

It was noticed in Section 3 of \cite{A1} that join semilattices with psudocomplemented upper sections admit a uniformly defined implication-like total operation which is a common extension of pseudocomplementations in all sections. Therefore, it can replace the family of section pseudocomplementations in any such a semilattice. The same idea was put in use in \cite{Ch-ext} for lattices. Another such operation for meet semilattices was defined in \cite{ChH}, but this definition appears implicitly already in  \cite[Section 2]{Sch}.

Our aim in this paper is to examine several ways of extending sectional pseudocomplementation on arbitrary posets. The next section contains the necessary general information on sectionally pseudocomplemented posets. Extended sectional pseudocomplementations (i.e., total extensions of sp-complementation), their basic properties and axiomatization is the subject of Section \ref{ext}. In Section \ref{related} we briefly review three related types of total operation on a poset ---  relative pseudocomplementation \cite{K2}, sectional j-pseudocomplemen\-tation \cite{C-impl}  and sectional pseudocomplementation in the sense of \cite{ChLP}, --- and show that generally these are not extended sp-complementations. 

In the rest of the paper we turn to concrete classes of posets with constructively described extended sp-comple\-men\-ta\-tion. Presented are, in particular, axiomatic characterizations of these classes and their subclasses formed by semilattices or lattices. The kind of extensions we deal with in Section \ref{elim} (the ``natural extensions'') is obtained  by adapting the technical definition of sp-complementation assumed in Section \ref{prelim}: we eliminate from the latter an implicit restriction that forces the defined operation to be partial. Then (in Section \ref{joinex}) we extend sp-complementation by means of a special rule which generalizes the one suggested in \cite{A1} for upper semilattices; it provides an adequate version of sectional j-pseudocomplementation. In Section \ref{another} we introduce a wide family of extended sp-complementations; it includes that   discussed in Section \ref{elim} and also an adequate variant of the one studied in \cite{ChLP}. Finally, in the concluding section we outline some proposals for further work.

\section{Preliminaries: sectional pseudocomplementation} \label{prelim}

In a poset $P := (P,\le)$, we denote by $[p)$ (by $(p]$) its upper (respectively, lower] section $\{x \in P\colon x \ge p\}$ (respectively $\{x \in P\colon x \le p\}$) determined by an element $p \in P$, and by $[p,q]$, the closed segment $\{x\colon p \le x \le q\}$; thus, $[p,q] = \varnothing$ when $p \not\le q$. Given a subset $X$ of $P$, we denote by $U(X)$ (by $L(X)$) the set of all upper (respectively, lower) bounds of $X$; thus, $U(\varnothing) = P = L(\varnothing)$. For elements $x,y,p \in P$, we write $x \perp_p y$ to mean that $x$ and $y$ are disjoint with respect to $p$, i.e., that there is no lower bound of $x$ and $y$ in $[p)$ distinct from $p$:
\[
x \perp_p y \tiff [p,x] \cap [p,y] \subseteq \{p\}.
\]
When really $p \le x,y$, this means that $p$ is the single, hence, the greatest lower bound of $x$ and $y$ in $[p)$ (in symbols, $p = x \wedge_p y$); if otherwise, the condition $x \perp_p y$ is trivially fulfilled. Thus,
\[
x \wedge_p y = z \tiff [p,x] \cap [p,y] = [p,z] \neq \varnothing.
\]
Likewise, the notations $x \vee y = z$ and $x \wedge y = z$ mean that the join  resp., the meet of $x$ and $y$  exists and equals to $z$. $P$ is an \emph{upper} (\emph{lower}) \emph{semilattice} if $x \vee y$ (resp., $x \wedge y$) exists for all $x,y \in P$. A \emph{join}, resp., \emph{meet semilattice} is then an algebra with an appropriate binary operation.

The pseudocomplement of $x$ in a section $[y)$, which we denote by $x * y$, can now be characterized as follows:
\begin{equation}    \label{x*y=}
x * y = \max \{u\colon u \wedge_y x  = y\}
\end{equation}
or, more explicitly, 
\begin{equation}    \label{=x*y=}
(x*y) \wedge_y x = y \mb{ and,}
\tforall u, \tif u \wedge_y x = y, \tthen u \le x*y.
\end{equation}
It follows from any of these characteristics that $x*y$ is defined only if $x \in [y)$, is unique in this case and belongs to $[y)$.

\begin{definition}    \label{sp-d}
A poset $P$ is said to be \emph{sectionally pseudocomplemented} (\emph{sp-comple\-mented}, or even an \emph{sp-poset} for short) if the pseudocomplement $x*y$ exists in it for all $x,y$ with $x \in [y)$. The partial operation $*$ arising on $P$ in this way is called a \emph{sectional pseudocomplementation}, or an \emph{sp-complementation} for short.
\deq
\end{definition}

There are two other useful characterizations of sectional pseudocomplements.

\begin{lemma}
Let $x,y \in P$. The  condition \eqref{=x*y=} is equivalent to any of the following ones:
\begin{gather}
 \label{<x*y}
x \in [y) \mb{ and,} \tforall u \in [y),\, (u \le x * y \tiff u \perp_y x), \\
\label{in*}
x \in [y) \mb{ and,} \tforall u,\, y \le u \le x*y \tiff u \wedge_y x = y.
\end{gather}
\end{lemma}
\begin{proof}
\eqref{=x*y=} $\im$ \eqref{in*}
Clearly, \eqref{=x*y=} implies the first clause of \eqref{in*} and the \textsl{if} part of the second one. On the other hand, if $y \le u \le x*y$, then, by \eqref{x*y=}, $u \le u'$ for some $u'$ such that $u' \wedge_y x = y$. As $y \le u$, it follows that $u \wedge_y x = y$. So \eqref{in*} follows from \eqref{=x*y=}.

\eqref{in*} $\eq$ \eqref{<x*y}
Recall that the equality $u \wedge_y x = y$ in \eqref{in*} is equivalent to the conjunction of the inequalities $y \le x$, $y \le u$ and the statement $u \perp_y x$. The inequality $y \le x$ may be omitted, as it follows from the first clause $x \in [y)$ of \eqref{in*}. Next, for any propositions $A,B,C$, the implication (if $A$, then ($B$ iff $C$)) is logically equivalent to the equivalence (($A$ and $B$) iff ($A$ and $C$)). Therefore \eqref{in*} is equivalent to \eqref{<x*y} (think the inequality $y \le u$, which means that $u \in (y]$, as $A$).

\eqref{in*} $\im$ \eqref{=x*y=}
Assume \eqref{in*}. From its second clause, $y \le x*y$ iff $y \le x$  (put there $u := y$). As $y \le x$ by the first clause, it follows that $y \le x*y$. Then, putting $u := x*y$ in the second clause, we come to the first clause in \eqref{=x*y=}. The second clause of \eqref{=x*y=} evidently follows from \eqref{in*}.
\end{proof}

\begin{lemma}    \label{sp-prop}
Sectional pseudocomplementation  has the following properties (supposing that the sectional pseudocomplements involved in each of them exist):
\begin{myenum}
\item   
$y \le x*y$,
\item   
$(x*y) \perp_y x$,
\item   
if $x \le y$, then $y*z \le x*z$,
\item   
if $x \le y*z$, then $y \le x*z$,
\item   
$x \le (x*y)*y$,
\item   
$y \le (x*y)*y$,
\item   
if $y \le x$, then $x \le y*y$,
\item   
if $y \le x$, then $(y*y)*x = x$,
\item   
$((x*y)*y)*y = x*y$,
\item   
if $x \le x*y$, then $x \le y$,
\item   
if $y \le x$, then $x*x = y*y$,
\item   
if $x*y = y*y$, then $x = y$,
\item   
if $z$ is a maximal lower bound of $x$ and $y$, then $x \le y*z$.
\end{myenum}
\end{lemma}
\begin{proof}
Assume that a poset $P$ is sectionally pseudocomplemented.

(a) Immediately from the first equality in \eqref{=x*y=}.

(b) By \eqref{=x*y=}, the first clause. 

(c) Suppose that  $x \le y$ and that $x*z$ and $y*z$ exist; so $x,y \in [z)$. As $(y*z) \perp_z y$ by (b), then also $(y*z) \perp_z x$. Now $y*z \le x*z$ in view of \eqref{<x*y} (and (a)).

(d) By \eqref{<x*y}, since $x \perp_z y$ iff $y \perp_z x$, when $x,y \in [z)$.

(e) By (d) (and (a)).

(f) An instance of (a).

(g) By (d) from (a), since $x*y$ is defined by the supposition.

(h) Suppose that $y \le x$; by (g), then $x \le y*y$ and the pseudocomplement $(y*y)*x$ is defined. Thus, $y \le x \le (y*y)*x$ by the supposition and (a), and then $(y*y)*x \le y*y$ by (g). Therefore, $(y*y)*x$ is a lower bound of both $(y*y)*x$ and $y*y$ in $[x)$. But  $((y*y)*x)  \perp_x (y*y)$ from (b), and we  conclude that $(y*y)*x = x$, as required.

(i) Suppose that $x*y$ is defined. Then $y \le x$. $y \le x*y$ by (a), and so $x*y \le ((x*y)*y)*y$ by (e) (and (a)). For the reverse inequality, apply (c) to (e).

(j) If $x \le x*y$, then $y \le x$ and $x$ is a lower bound of $x$ and $x*y$ in $[y)$. But $x*y \perp_y x$ by (b); so, $x = y$.

(k) By (a), $x \le x*x$; so, if $y \le x$, then $y \le x*x$. Using (g), $x*x \le y*y$, whence $x \le y*y$. Applying (g) once more, $y*y \le x*x$.

(l) Suppose that $x*y = y*y$ and, therefore, $y \le x$. Then $x \le y*y$ by (g), so that $x \le x*y$ and, by (j), $x \le y$.

(m) An element $z$ is a maximal lower bound of $x$ and $y$ if and only if
($\alpha$) $z \le x,y$ and ($\beta$) for all $v \le x,y$, $v = z$ whenever $z \le v$. But ($\beta$) means that $x \perp_z y$. Thus ($\alpha$) and ($\beta$) imply that $x \le y*z$; see \eqref{<x*y}.
\end{proof}

So, the conjunction of ($\alpha$) and ($\beta$) in the the proof of (m)  means that $x \wedge_z y = z$. Therefore,
\begin{equation}    \label{mlb}
p \text{ is a maximal lower bound of } x \tand y\ \tiff\ x \wedge_p y = p.
\end{equation}
The following observations on the structure of an sp-poset are consequences of Lemma \ref{sp-prop}(a,g).

\begin{corollary}    \label{p*p}
For every $p$ in a sectionally pseudocomplemented poset, $p*p$ is its maximal element, which is also greatest in $[p)$. Different maximal elements of the poset do not have any common lower bound. The poset has a greatest element if and only if it is up-directed.
\end{corollary}
\begin{proof}
Let $P$ be an sp-poset. If $p*p \le x$, then $p \le x$ by Lemma \ref{sp-prop}(a), and if $p \le x$, then  $x \le p*p$  by Lemma \ref{sp-prop}(g); this justifies the first statement. For the next one, notice that if  elements $a,b$ are maximal and $c$ is their lower bound, then $a = c*c = b$ in view of Lemma \ref{sp-prop}(g). Finally, if $P$ is up-directed, then there is only one maximal element in it, which is then the greatest one; the `only if'  part of the third statement is self-evident. 
\end{proof}

A poset is said to be \emph{sectionally bounded} if every section $[p)$ in it is bounded from above and is therefore a closed segment. We denote the greatest element of the section [p) by $1_p$; thus, $[p) = [p,1_p]$. So, every sp-poset is sectionally bounded, and the following holds in it:
\begin{equation}   \label{1p}
1_p = p*p .
\end{equation}
Furthermore, if a sectionally pseudocomplemented poset is an upper or lower semilattice, then it necessarily has a greatest element.

We end the section with a simplified axiomatic description of posets equipped with a sectional pseudocomplementation. They are considered as partial algebras of kind $(P,*)$, where $P = (P,\le)$ is a poset  and $*$ is a partial binary operation defined for a pair $(x,y)$ if and only if $y \le x$; let us call the ordered algebras of this kind \emph{star posets} or, in short, \emph{$*$-posets}. When a star poset has the greatest element 1, we shall present it, if necessary, as a triple $(P,*,1)$.
If the operation $*$ satisfies, for every $p$, also the condition
``if $x \in [p)$, then $x*p \in [p)$'',
then every section $[p)$ is closed under $*$ in the sense that $x*y \in [p)$ whenever $x,y \in [p)$ and $y \le x$. Such an operation may be called \emph{sectional}. For instance, sp-complementation is a sectional operation---recall Lemma \ref{sp-prop}(a).

\begin{definition}    \label{$*$--sp}
A \emph{poset with sectional pseudocomplementation} is a $*$-poset $(P,*)$ where $*$ is an sp-complementation on $P$. We usually shall use a shorter term \emph{SP-poset} for such a $*$-poset, and denote the class of these structures by \SP.
\deq
\end{definition}

\begin{example}    \label{sp-e}
Let $P$ be a poset depicted below. It is sectionally pseudocomplemented: the partial operation $*$ on it defined by the table at the middle:
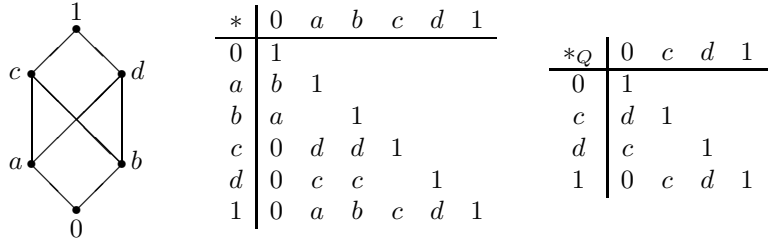
\setlength\unitlength{0.6mm}
\[
\hspace*{1,4cm}
\begin{picture}(10,20)(25,19)
\put(10,40){\circle*{2}}
\put(0,30){\circle*{2}} \put(20,30){\circle*{2}}
\put(0,10){\circle*{2}} \put(20,10){\circle*{2}}
\put(10,0){\circle*{2}}
\put(0,30){\line(1,1){10}}
\put(20,30){\line(-1,1){10}}
\put(0,30){\line(1,-1){20}}
\put(20,30){\line(-1,-1){10}}
\put(0,10){\line(1,1){20}}
\put(20,10){\line(-1,1){20}}
\put(0,10){\line(1,-1){10}}
\put(20,10){\line(-1,-1){10}}
\put(0,10){\line(0,1){20}}
\put(20,10){\line(0,1){20}}
\put(8.5,42){1}
\put(-5,29){$c$} \put(22,29){$d$}
\put(-5,9){$a$} \put(22,9){$b$}
\put(8.5,-6){\mbox{0}}
\end{picture}
\quad
\begin{array}{c|cccccc}
* &0&a&b&c&d&1 \\
\hline
0 &1& & & & & \\
a &b&1& & & & \\
b &a& &1& & & \\
c &0&d&d&1& & \\
d &0&c&c& &1& \\
1 &0&a&b&c&d&1
\end{array}
\qquad
\begin{array}{c|cccc}
*_Q &0&c&d&1 \\
\hline
0 &1& & & \\
c &d&1& & \\
d &c& &1& \\
1 &0&c&d&1
\end{array}
\]
is an sp-complementation. Its subset $Q := \{0,c,d,1\}$ is both sp-complemented (see the table for $*_Q$) and closed under $*$. So, $(P,*)$ is an SP-poset, but $(Q,*)$ is not.
\deq
\end{example}

The following characterization of SP-posets is a modified version of Lemma 7 in \cite{C-annih},  which was given there without proof.

\begin{theorem}    \label{sp-char}
The class \SP\ is characterized by the axioms  

\begin{myenum}
\thitem{sp1}
if $x \le y$, then $y*z \le x*z$,
\thitem{sp2}
if $x \le x*y$, then $x \le y$,
\thitem{sp3}
if $z$ is a maximal lower bound of $x$ and $y$, then $x \le y * z$,
\end{myenum}
where the elements $x*y, x*z, y*z$  involved in the first two of them are supposed there to exist.
\end{theorem}
\begin{proof}
Let $(P,*)$ be a $*$-poset. Then the statements (sp1)--(sp3) are fulfilled, see Lemma \ref{sp-prop}(c,j,m). Now assume that they, conversely, are fulfilled in $P$.

First, then also the inequality
\begin{enumerate}
\thitem{*} $x \le (x*y)*y$
\end{enumerate}
from Lemma \ref{sp-prop}(e) holds. Indeed, suppose that $x*y$ is defined and therefore $y \le x$. By (sp3), then $y \le x*y$; so, $(x*y)*y$ is defined and $y$ is a lower bound of $x$ and $x*y$. If $v$ is another lower bound and $y \le v$, then $v \le x*y \le v*y$ by (sp1) and further $v \le y$ by (sp2). Thus, $y$ is even a maximal lower bound of $x$ and $x*y$, whence $x \le (x*y)*y$ by (sp3).

Next, consider again any element $x*y$; for it, the condition \eqref{<x*y} should be fulfilled. Of course, $x \in [y)$. Suppose further that $u \in [y)$. If $u \perp_y x$, then $y$ is a maximal lower bound of $u$ and $x$, hence, $u \le x *y$ by (sp3). If, conversely, $u \le x*y$, then $u \perp_y x$ follows: if $y \le v \le u,x$, then $v \le x, x*y$ (immediately), $(x*y)*y \le v*y$ by (sp1), $v \le x \le v*y$ due to (*), and finally $v = y$ by (sp2). So,  \eqref{<x*y} holds, and $(P,*)$ is an SP-poset.
\end{proof}

Thus, the conjunction of the axioms (sp1) and (sp2) is equivalent, in the presence of (sp3), to the first clause in \eqref{=x*y=}  (see \eqref{mlb}).

\section{Extended SP-posets}  \label{ext}

\subsection{Extended sp-complementation.} 
Ordered algebras $(P,\to)$, where $P$ is a poset and $\to$ is a total binary operation on $P$, will be called \emph{arrow posets} (we thus adapt the notion of arrow semilattice from  \cite{C-wrpseudo}). When such an algebra has a greatest element 1 and it is necessary to stress this, we present the arrow poset as a triple $(P,\to,1)$.

\begin{definition}    \label{esp-e}
An \emph{extended sectional pseudocomplementation} (\emph{esp-complementa\-tion}) is a binary  operation $\to$ on a poset $P$ such that $x \to y$ is a pseudocomplement of $x$ in $[y)$ for all $x$ and $y$ with $y \le x$. An arrow poset $(P,\to)$, where $\to$ is an esp-complementation, is called a \emph{poset with esp-complementation}, or an \emph{extended SP-poset} for short, also an \emph{ESP-poset}. It is an \emph{extension of an SP-poset $(P,*)$} if $\to$ is an extension of $*$. The class of all ESP-posets is denoted by \ESP.
\deq
\end{definition}

Thus, every SP-poset can be extended to an ESP-poset, usually in various ways. A remark: in this paper, we do not speak of  an extension of a poset as of a poset in which the former one is embedded.

Theorem \ref{sp-char} immediately leads to the following axiomatic characterization of  ESP-posets.

\begin{proposition}    \label{esp-char}
The class \ESP\ is characterized by axioms \\[2pt]
\textup{(esp1)}
if $z \le x \le y$, then $y \to z \le x \to z$, \\
\textup{(esp2)}
if $y \le x \le x \to y$, then $x \le y$, \\
\textup{(esp3)}
if $z$ is a maximal lower bound of $x$ and $y$, then $x \le y \to z$.
\end{proposition}
\begin{proof}
Let $P$ be a poset, and let $\to$ be an operation on $P$ satisfying the three axioms. The axiom (esp3) implies that $y \le x \to y$ whenever $y \le x$. Then the restriction $*$ of $\to$, characterized by
\[
x*y = z\ \tiff \ x \to y = z \tand y \le x,
\]
is a sectional operation on $P$. Thus, $*$ is, by Theorem \ref{sp-char}, an sp-complementation, and therefore $\to$ itself, an esp-complementation.

If, conversely, $\to$ is an extension of an sp-complementation $*$, then $x \to y = x*y$ whenever $x \le y$. Consequently, the axioms (esp1)--(esp3) reflect the conditions (sp1)--(sp3) and therefore are fulfilled in $(P,\to)$.
\end{proof}

This result implies that a poset is sp-complemented if and only if it admits an operation $\to$ satisfying (esp1)--(esp3). For further references, we now list some properties of extended sp-complemen\-tation.

\begin{lemma}    \label{esp-prop}
In an ESP-poset $(P,\to)$,
\begin{myenum}
\item   
if $y \le x$, then $y \le x \to y$,
\item   
if $y \le x$, then $x \le (x \to y) \to y$,
\item   
if $z \le x,y$ and $x \le y \to z$, then $y \le x \to z$,
\item   
if $y \le x$, then $y \le (x \to y) \to y$,
\item    
if $y \le x$, then $((x \to y) \to y) \to y = x \to y$,
\item   
$x \to x = 1_x$.
\item   
if $y \le x$, then $x \le 1_y$,
\item   
if $y \le x$, then $1_y \to x = x$,
\item   
if $y \le x$, then $1_x$ = $1_y$.
\end{myenum}
\end{lemma}
\begin{proof}
Items (a)--(e) and (g)--(i) follow from Lemma \ref{sp-prop}(a,e,d,f,i,g,h,k), respectively, and (f) follows from \eqref{1p}. 
\end{proof}

\subsection{Implicative ESP-posets.}
Of particular interest are those ESP-posets where the order $\le$ is definable in terms of $\to$: they can be presented as ordinary (i.e., unordered) binary algebras.
Implicative algebras of  \cite {Ra} can be characterized as algebras $(A,\to,1)$ where the relation $\le$ introduced by the condition ``$x \le y$ iff $x \to y = 1$'' is a partial order with $1$ the greatest element. The arrow poset arising in this way from such an algebra can also be called implicative. In the definition below, we extend this latter attribute to  certain arrow posets which may have not the greatest element.

\begin{definition}
An arrow poset $(P,\to)$ is said to be \emph{left implicative} (\emph{right implicative}) if it is sectionally bounded and, for all $x,y$,
\begin{equation}    \label{le/to}
x \le y \tiff x \to y = 1_x,
\end{equation}
respectively,
\begin{equation}    \label{le/to2}
x \le y \tiff x \to y = 1_y,
\end{equation}
and \emph{implicative} if it is both left and right implicative.
An extension of a sectional operation on a poset  is left implicative, right implicative or {implicative} if so is the corresponding arrow poset.
\deq
\end{definition}

Thus, in a left-implicative (right-implicattive) ESP-poset,
\[
\mbox{if $x \to y = 1_x$, respectively,  $x \to y = 1_y$, then $1_x = 1_y$};
\]
see Lemma \ref{esp-prop}(i). Of course, an SP-poset generally may have extensions that are neither left nor right implicative. Generally, none of the conditions \eqref{le/to} and \eqref{le/to2} implies the other one.

\begin{example}    \label{impl}
The poset $P$ consisting of two maximal chains $a < b$ and $c < d$ is sectionally pseudocomplemented; see the first table. So, $b = 1_a = 1_b$ and $d = 1_c = 1_d$. The  operation $\to_1$ is a left-implicative extension of $*$\,; however, it is not right 
\[
\begin{array}{c|ccccc}
* &a&b&c&d \\
\hline
a   &b& & &  \\
b   &a&b& &  \\
c   & & &d& \\
d   & & &c&d
\end{array}
\quad
\begin{array}{c|ccccc}
\to_1 &a&b&c&d \\
\hline
a   &b&b&d&d \\
b   &a&b&d&d \\
c   &b&b&d&d \\
d   &b&b&c&d
\end{array}
\quad 
\begin{array}{c|ccccc}
\to_2 &a&b&c&d \\
\hline
a   &b&b&b&b \\
b   &a&b&b&b \\
c   &d&d&d&d \\
d   &d&d&c&d
\end{array}
\quad
\begin{array}{c|ccccc}
\to_3 &a&b&c&d \\
\hline
a   &b&b&a&c \\
b   &a&b&c&a \\
c   &c&a&d&d \\
d   &a&c&c&d
\end{array}
\]
implicative, as $c \to a = 1_a$ though $c \nleq a$. Likewise, the operation $\to_2$ is a right-implicative extension of $*$, but is not left implicative, as $c \to a = 1_c$. The operation $\to_3$ is an implicative extension of $*$.
\deq
\end{example}

\begin{proposition}
A right-implicative ESP-poset in which $y \le x \to y$ for all $x$ and $y$ is implicative. 
\end{proposition}
\begin{proof}Let $(P,\to)$ be an ESP-poset satisfying both presuppositions. If $x \le y$ then $x \to y = 1_y = 1_x$ by Lemma \ref{esp-prop}(i). If $x \to y = 1_x$, then $y \le 1_x$, whence $1_x \le 1_y$ by Lemma \ref{esp-prop}(g) and, further, $1_x = 1_y$. So, $x \to y = 1_y$, whence $x \le y$. Therefore, $(P,\to)$ is also left implicative.
\end{proof}

 The next theorem shows that there is an abundance of implicative ESP-posets.

\begin{theorem}
Every SP-poset admits an implicative extension.
\end{theorem}
\begin{proof}
Let $(P,*)$ be an SP-poset. Recall that $p*p$ is the greatest element $1_p$ in its section $[p)$.
Consider a total binary operation $\to$ on $P$ is defined as follows:
\begin{equation}    \label{pure}
x \to y := \left\{
\begin{array}{l}
x * y \mbox{ if } y \le x, \\
1_x \mbox{ if } x < y,    \\
y \mbox{ otherwise.}  \\
\end{array}
\right.
\end{equation}
The arrow poset $(P,\to)$ is a left-implicative extension of $(P,*)$. Indeed, $\to$ is an extension of $*$, and the \textsl{only if} part of (\ref{le/to}) also is an immediate consequence of (\ref{pure}) and \eqref{1p}. Furthermore, suppose that $x \to y = 1_x$; to prove that $x \le y$ it suffices to consider the case when $y \neq 1_x$. Then $x \to y \neq y$, and (\ref{pure}) implies that $y$ must be comparable with $x$. Now, if $y < x$, then $x*y = x \to y = 1_x = 1_y$ (by \eqref{pure}, the supposition and Lemma \ref{esp-prop}(i)), whence $x = y$ (by Lemma \ref{sp-prop}(l) and \eqref{1p}), a contradiction. Thus $x \le y$, and the \textsl{if} part of (\ref{le/to}) is also fulfilled.

By Lemma \ref{esp-prop}(i), the value $1_x$ on the second line of the definition (\ref{pure}) can equivalently be replaced by $1_y$. A similar proof then shows that the ESP-poset $(P,\to)$ is also right implicative.
\end{proof}

\begin{definition}
We call an implicative extension of an SP-poset \emph{pure} if it is  defined via the rule (\ref{pure}).
\deq
\end{definition}

\begin{example}    \label{sp-e1}
Let $(P,*,1)$ be the sp-poset from Example \ref{sp-e}. The arrow poset $(P,\to,1)$ with the table
\[
\begin{array}{c|cccccc}
\to&0&a&b&c&d&1 \\
\hline
0 &1&1&1&1&1&1 \\
a &b&1&b&1&1&1 \\
b &a&a&1&1&1&1 \\
c &0&d&d&1&d&1 \\
d &0&c&c&c&1&1 \\
1 &0&a&b&c&d&1
\end{array}
\]
is its pure implicative extension. Clearly, the implicative extension  $(P,\to_1,1)$ described in the preceding example is not pure.
\deq
\end{example}

\begin{remark}  \label{notQV}
Any  (left-, right-)implicative ESP-poset may be considered also as an ordinary algebra where  \eqref{le/to}, respectively, \eqref{le/to2} defines a derived order relation; let us call such an algebra an (\emph{left-, right-})\emph{implicative ESP-algebra}. The classes of all such algebras are not equational; even more, Example \ref{sp-e1} presents an implicative ESP-algebra whose subalgebra determined by the subset $Q := \{0,c,d,1\}$ does not belong to this class; see Example \ref{sp-e}. Being therefore not closed under forming subalgebras, the class of (left-, right-)implicative ESP-algebras cannot be axiomatized by quasi-equations (alias equational implications) and therefore is not a quasivariety (an implicational class); see \cite[\S63, Theorem 3]{UA}. 
\deq
\end{remark}

\subsection{ESP-semilattices}
The following proposition shows that in semilattices the axiom (esp3) can be improved.

\begin{proposition}    \label{glb}
The following conditions on $x,y,z$ are equivalent in any upper or lower semilattice:
\begin{myenum}
\item $z$ is a maximal upper bound of $x$ and $y$,
\item $z$ is a greatest lower bound of $x$ and $y$,
\item  $z$ is a greatest lower bound of $x$ and $y$ in $[z)$.
\end{myenum}
\end{proposition}
\begin{proof}
Equivalence of (a) and (c) (even in an arbitrary poset) is stated in \eqref{mlb}, while equivalence of (a) and (b) in upper semilattioces is evident. Now assume that $P$ is a lower semilattice, and assume (a). So, $z \le x,y$, whence $z \le x \wedge y$. But $x \wedge y$ is a lower bound of $x$ and $y$; so $z = x \wedge y$ by maximality of $z$. Thus (a) implies (b); the converse is trivial.
\end{proof}

Notice that in the poset from Example \ref{sp-e}, which is not a semilattice, the meet  $c \wedge d$ does not exist, whereas $c \wedge_b d = b$.
 
For meet semilattices, we can obtain even more than one improved axiom. As usual, an inequality $a \le b$ in a meet semilattice can be treated as a shorthand for the equality $a \wedge b = a$ when necessary.

\begin{theorem} \label{ESP^}
The class \ESP$^\wedge$ of meet semilattices with esp-comple\-mentation is a variety determined by the semilattice axioms and the identities
\begin{enumerate}
\thitem{esp$^\wedge$1}  $x \wedge (x \to (x \wedge y)) = x \wedge y$, 
\thitem{esp$^\wedge$2}  $x \le y \to (x \wedge y)$.
\end{enumerate}
\end{theorem}

\begin{proof}
Let $(P,\wedge,\to)$ be a arrow meet semilattice. If it is an algebra from \ESP$^\wedge$, then, for $y \le x$, the first item of \eqref{=x*y=} yields  the equality $(x \to y) \wedge_y x = y$. Using  Proposition \ref{glb}, then
\[
\tif y \le x, \tthen x \wedge (x \to y) = y;
\]
this conclusion evidently implies (esp$^\wedge$1). In its turn, (esp$^\wedge$2)
follows from (esp3). Conversely, assume that (esp$^\wedge$1) and (esp$^\wedge$2) hold in $(P,\wedge, \to)$. Then (esp2) easily follows from (esp$^\wedge$1), and (esp3) follows from (esp$^\wedge$2) by Proposition \ref{glb}. To prove (esp1), suppose that $z \le x \le y$. Then $z$ is a lower bound of $x$ and (by (esp$^\wedge$2)) $y \to z$. It is even their greatest lower bound: if $v \le x, y \to z$, then $v \le y$, $y \to (y \wedge z)$ (see the supposition), and (esp$^\wedge$1) implies that $v \le z$. So, $y \to z \le x \to z$ by (esp$^\wedge$2).
\end{proof}

\begin{corollary} \label{ESP^v}
The class \ESP$^{\wedge\vee}$ of lattices with esp-complementation is a variety determined by the lattice axioms and the axioms \textup{(esp$^\wedge$1), (esp$^\wedge$2)}.
\end{corollary}

\begin{remark}  \label{sp=wr}
Evidently, the axiom  (esp$^\wedge$1)  is equivalent to the weaker identity 
\[
\mbox{\upshape{(esp$^\wedge$1')} $x \wedge (x \to (x \wedge y)) \le y$.}
\]
According to \cite[Proposition 7]{C-wrpseudo}, axioms (esp$^\wedge$1') and (esp$^\wedge$2) determine the class \EWR$^\wedge$ of meet semilattices with extended weak relative pseudocomplementation. Thus \ESP$^\wedge$ = \EWR$^\wedge$ (this essentially follows also from Proposition 2 in \cite{C-wrpseudo}). In particular, the congruence and filter properties of meet \EWR-semilattices discovered in \cite{C-wrpseudo} are possessed also by semilattices in \ESP$^\wedge$. For more on weak relative pseudocomplementation, see Section \ref{related}.3 below. 
\deq
\end{remark}

\section{Some related kinds of extended sectional operations}   \label {related}

In this section, we briefly review three classes on arrow posets known from the literature and related by name to sectional pseudocomplementation. We show that actually none of them is a subclass of \ESP.

\subsection{Relative pseudocomplementation \cite{K2}}

Relative  pseudocomplementation (\emph{rp-complementation} in short) is defined as follows:
\[
x \to y := \max\{u\colon (u] \cap (x] \subseteq (y]\}.
\]
As $(y] \cap (x] \subseteq (y]$, evidently, $x \to y \in [y)$ for $x \in [y)$; so $(P,\to)$ is an extension of a unique $*$-poset $(P,*)$, in which $*$ is a sectional operation characterized by
\[
\tforall u,\; u \le x * y \tiff (u] \cap (x] \subseteq (y] .
\]
\begin{example}    \label{rpc-e}
Example 4 in \cite{ChL} presents a poset with relative pseudocomplementation $(P,\to,1)$, where  $P$ is the poset  from Example \ref{sp-e}  above and $\to$ has the operation table given below. The  other table shows the sectional operation $*$ which
\newline\vspace*{-1em}
\[
\begin{array}{c|cccccc}
\to &0&a&b&c&d&1 \\
\hline
0 &1&1&1&1&1&1 \\
a &b&1&b&1&1&1 \\
b &a&a&1&1&1&1 \\
c &0&a&b&1&d&1 \\
d &0&a&b&c&1&1 \\
1 &0&a&b&c&d&1
\end{array}
\qquad
\begin{array}{c|cccccc}
* &0&a&b&c&d&1 \\
\hline
0 &1& & & & &  \\
a &b&1& & & &  \\
b &a& &1& & &  \\
c &0&a&b&1& &  \\
d &0&a&b& &1&  \\
1 &0&a&b&c&d&1
\end{array}
\]
is the restriction of $\to$. A comparison of this table with that of Example \ref{sp-e} shows that $*$ here is not an sp-complementation. Therefore, on one side, $\to$ is not an esp-complementation; on the other side, the sp-complementation of Example \ref{sp-e} is not a restricted rp-complementation.
\deq
\end{example}

Thus the class of all posets with rp-complementation is not a subclass of \ESP. Moreover, as every relatively pseudocomplemented poset is distributive, it follows that the sp-complementation on a distributive poset not necessarily can be extended to an rp-complementation---while it can on a lower semilattice that is only modular. See \cite[Section 2]{C-wrpseudo} for details.

\subsection{Sectional j-pseudocomplementation \cite{C-impl}}
\label{sjc-posets}

A certain class of upper semilattices with extended sectional pseudocomplementation and unit  was characterized in \cite[Theorem 1]{C-impl} by the axioms
\begin{gather*}
\mbox{(j1) if $x \le y \to z$, then $y \le x \to z$, \qquad
(j2) if $x \le x \to y$, then $x \le y$,} \\
\mbox{(j3) if the meet of $x$ and $y$ exists, then $x \le y \to (x \wedge y)$.}
\end{gather*}
Arbitrary arrow posets with 1 satisfying these axioms were then said to be \emph{sectionally j-pseudocom\-plemented} and considered as particular extended SP-posets. However, a crucial point in the proof of the theorem was the fact that the subcondition $u \perp_y x$ appearing in the definition of sectional pseudocomplementation used there (essentially the condition \eqref{<x*y} above) is equivalent, in an upper semilattice, to the equality $u \wedge x = y$ when $y \le x,u$. This means that the axioms (j1)--(j3) may be not adequate for posets that are not  upper semilattices (cf.\ the note following Proposition \ref{glb}). And indeed, there are sectionally j-pseudocomplemented posets which actually are not ESP-posets.

\begin{example}    \label{sjpc-e}
The arrow poset from Example \ref{rpc-e}, being relatively pseudocomplemented, satisfies the axioms (j1)--(j3) and is therefore sectionally j-pseudo\-comple\-mented. However, it is not an extension of an SP-poset as we already know.
\qed
\end{example}

An improved system of axioms excluding such examples will be presented in Section \ref{joinex}.

\subsection{Sectional pseudocomplemention as a total operation \cite{ChLP}}

We now turn to an approach to sectionally pseudocomplemented posets which was motivated by the specific definition of sectionally pseudocomplemented lattices mentioned in Introduction. In \cite[Definition 2.1]{ChLP}, a poset $P$ was said to be sectionally pseudocomplemented if, for all $x,y$, there exists a greatest $u \in P$ such that
\[  \tag{$\ast$}
L([x) \cap [y)) \cap (u] = (y] .
\]
This greatest element was termed a \emph{sectional pseudocomplement of $x$ with respect to $y$}. We denote it here by $x \to y$; the corresponding arrow poset $(P,\to)$ also was said to be sectionally pseudocomplemented. To avoid a terminological conflict, we say that it is  a \emph{CLP-pose} ('C', 'L', 'P' for the authors'   names). Thus, the following is an explicit definition of the operation $\to$:
\begin{equation}    \label{CLP}
x \to y = \max\{u\colon L([x) \cap [y)) \cap (u] = (y]\}.
\end{equation}
For instance, all relatively pseudocomplemented posets are  CLP-sec\-tionally pseudocomplemented (\cite[Theorem 2.8]{ChLP}). 

Let us make a closer look at the definition \eqref{CLP}. If $y \le x$, then the above defining equality ($\ast$) reduces to
$(x] \cap (u] = (y]$. Thus, writing $x * y$ for $x \to y$ in this case,
\[
x * y = \max\{u\colon (u] \cap (x] = (y]\}.
\]
This equality differs from that presenting sp-complementation, \eqref{x*y=}. It actually defines a weak relative pseudocomplementation (\emph{wrp-com\-plementation}) in the sense of \cite[Section 3]{C-annih}. Both kinds of local pseudocomplementation  coincide in meet  semilattices (see Remark \ref{sp=wr} above in Section \ref{ext}) and upper semilattices (due to Proposition \ref{glb}), but things again are not so good in arbitrary posets: there are posets where the wrp-complementation is not an sp-complementation (and conversely). This is immediate from Example \ref{rpc-e}, as the operation $*$ there proves to be a wrp-complementation. (Even more, every rp-complementation is an extended wrp-complementation: when $y \le x$,
$\max\{u\colon (u] \cap (x] \subseteq (y]\}
= \max\{u \ge y\colon  (u] \cap (x] \subseteq (y]\} = \max\{u\colon (u] \cap (x] = (y]\}$.)

So, the SP-poset from Example \ref{sp-e} is an instance of a bounded SP-poset that cannot be extended to an CLP-poset (recall Example \ref{rpc-e} and the note following \eqref{CLP}). What is of more significance, the arrow poset from Example \ref{rpc-e} is an instance of a CLP-poset that actually has sections which are not sp-complemented.

 The above notes show that a CLP-poset always is an ESP-poset if it is a lower semilattice (this conclusion repeats Proposition 2.3 of \cite{ChLP}). However, CLP-posets generally should be considered as extensions of posets with wrp-complementation rather than sp-complementation. A modification of the definition \eqref{CLP} that isolates only ESP-posets will be presented in Section \ref{another}.4.

\section{Extension of sp-complementation by relaxing its definition}
\label{elim}

According to \eqref{<x*y}, the pseudocomplement of $x$ in $(y]$ is, for $x \in [y)$, the greatest element of $[y)$ disjoint with $x$ in $(y]$.  A natural way to obtain a totally defined extension of sectional pseudocomplementation is to eliminate from  this description the restriction $x \in [y)$. 
Thus, consider the following description of an operation $\to$ on a poset $P$:
\begin{equation} \label{xtoy=}
x \to y = \max\{u \ge y\colon u \perp_y x\}.
\end{equation}
Equivalently,  
\begin{multline}    \label{=xtoy=}
\mb{(i)} y \le x \to y, \ \mb{(ii)} (x \to y) \perp_y x, \\
\mb{(iii)}  \tif u \ge y \tand u \perp_y x, \tthen u \le x \to y.
\end{multline}
One more useful characteristic of $\to$ is the equivalence
\[ 
\tforall u \in [y),\, u \le x \to y \tiff u \perp_y x,
\] 
which is comparable with the second clause in \eqref{<x*y}. 
Therefore, if the operation $\to$ in \eqref{xtoy=} turns out to be total, it is an extended sp-complementation on $P$. Notice that the clause (i) in \eqref{=xtoy=} is in fact superfluous: as always $y \perp_y x$, it follows from (iii) by substituting $y$ for $u$.

\begin{definition}
An operation $\to$ on a poset $P$ is called a \emph{natural} esp-complemen\-tation if it satisfies \eqref{xtoy=}; we also say in this case that the operation is a natural extension of the corresponding sp-complementation. The arrow poset $(P,\to)$ is then called a \emph{natural ESP-poset}, also a \emph{naturally extended} SP-poset, 
\deq
\end{definition}

In view of \eqref{=xtoy=}(i), the  element $x \to y$ in a natural ESP-poset $(P,\to)$ may be called a \emph{natural esp-complement of $x$ in $[y)$} (by abuse of language, it could even be called a {natural pseudocomplement of $x$ in $[y)$}).

The next theorem shows, in particular, that every SP-poset has a natural extension.

\begin{theorem} \label{nat1y}
An operation $\to$ on a poset $P$ satisfies \eqref{xtoy=} if and only if it is an extended sp-complementation and, for all $x,y$,
\begin{equation}    \label{!1}
x \to y = 1_y \mb{whenever} y \nleq x.
\end{equation}
\end{theorem}
\begin{proof}
Assume that the operation $\to$ satisfies \eqref{xtoy=}. As noted above, it is then an extension of an sp-complementation. Suppose further that $y \nleq x$ and therefore $[y,x] = \varnothing$; then $u \perp_y x$ for every $u \ge y$. Consequently, the equality \eqref{xtoy=} reduces, under this supposition, to $x \to y = \max\{u\colon u \ge y\} = 1_y$.

Now assume that, conversely, $\to$  is an extended sp-complementation and \eqref{!1} holds. Then $\to$ satisfies \eqref{xtoy=}: if $y \le x$, then $x \to y = x*y$, and \eqref{xtoy=} follows from \eqref{x*y=}, while  if $y \nleq x$, then $y \to x = 1_y = \max\{u \ge y\colon [y,u] \cap [y,x] \subseteq \{y\}\}$.
\end{proof}

Put in another way, the theorem says that the natural extension of an sp-complementation $*$ is given by the following rule:
\begin{equation}	\label{purerule}
x \to y := \left\{
\begin{array}{l}
x * y \mbox{ if } y \le x, \\
1_y \mbox{ otherwise}.
\end{array}
\right.
\end{equation} 
The class of natural ESP-posets can be characterized also by a set of axioms simpler by those provided by the theorem.

\begin{theorem}	\label{nat}
An operation $\to$ satisfies \eqref{xtoy=} if and only if it satisfies
the axioms
\begin{enumerate}
\thitem{nat1} if $x \le y$, then $y \to z \le x \to z$,
\thitem{nat2} if $y \le x \le x \to y$, then $x \le y$,
\thitem{nat3} if $z \le x$ and $x \perp_z y$, then $x \le y \to z$.
\end{enumerate}
\end{theorem}
\begin{proof}
Notice that the axiom (nat3) is a variant of the clause  (iii) of \eqref{=xtoy=}. Thus, it remains to prove that the conjunction of (nat1) and (nat2) is equivalent to the clause (ii) in the presence of (nat3) or (iii).

Now let $(P,\to)$ be an arrow poset. Assume that it satisfies the three axioms. Then the subcondition (ii) is fulfilled: if $y \le v \le x,x \to y$, then $v \le x \to y \le v \to y$ by (nat1) and, further, $v \le y$  by (nat2). So, \eqref{=xtoy=} holds. Conversely, assume that the three conditions in \eqref{=xtoy=} are fulfilled. By \eqref{=xtoy=}(ii), the inequalities $y \le x \le x \to y$ imply that $x = y$, i.e., (nat2) is fulfilled. To prove (nat1), suppose that $x \le y$ for some $x,y$ and choose any $z$. From (i), $z \le y \to z$. Next, $(y \to z) \perp_z x$: if $w \in [z, y \to z] \cap [z,x]$, then $w \in [z,y \to z] \cap [z,y]$  by the last supposition and $w = z$ by (ii). By (iii), finally  $y \to z \le x \to z$, as needed. So, all axioms are satisfied.
\end{proof}

Informative is also the following proof of the same theorem, which rests on Theorem \ref{nat1y}. For instance, it shows how the condition \eqref{!1} springs from (nat3).

\begin{proof}
Let $(P,\to)$ be an arrow poset.
First notice that (nat2) coincides with (esp2). Next, (nat1) is a conjunction of (esp1) and a condition 
\[	\tag{*}
\mbox{if $z \nleq x \le y$, then $y \to z \le x \to z$.}
\]
 Finally, the axiom (nat3) is equivalent to the conjunction of (esp3) and  \eqref{!1}. Indeed, it is 
a conjunction of statements
\par -- if $z \le y$, $z \le x$ and $x \perp_z y$  (i.e., if $x \wedge_z y$), then $x \le y \to z$, and 
\par -- if $z \nleq y$, $z \le x$ and $x \perp_z y$, then $x \le y \to z$.
\\  By \eqref{mlb}, the first one is equivalent  to (esp3). As $x \perp_z y$ trivially holds when $z \nleq y$, the second one actually states that if $z \nleq y$ and $z \le x$, then $x \le y \to z$; this implication is evidently equivalent to
\[	
\mbox{if $z \nleq y$, then $1_z \le y \to z$},
\] 
but implies also that if $z \nleq y$, then $y \to z \in [z)$ and, further, $y \to z \le 1_z$. Thus, the second displayed statement is equivalent to  \eqref{!1}. 

Now, if $(P,\to)$ satisfies the axioms (nat1)--(nat3), then it satisfies also (esp1)--(esp3) and  \eqref{!1}, i.e., is a natural ESP-poset. If, conversely, $(P,\to)$ is a natural ESP-poset, then (nat2) and (nat3) hold in $(P,\to)$ by the previous work. Moreover, then $y \to z \in [z)$  by \eqref{=xtoy=}(i), whence $y \to z \le 1_z$. It follows by \eqref{!1} that (*) is fulfilled, and therefore (nat1) also holds. 
\end{proof}

It turns out that a natural ESP-poset is rarely left or right implicative. The arrow poset $(P,\to_1)$ in Example \ref{impl} illustrate the following observations.

\begin{proposition}
A natural ESP-poset is 
\begin{myenum} 
\item left implicative if and only if every its lower section is a chain,
\item right implicative if and only if it is a chain.
\end{myenum}
\end{proposition}
\begin{proof}
Let $(P,\to)$ be a natural ESP-poset. 

(a) First assume that $P$ fulfills the condition (then it is a chain or an union of disjoint chains), and suppose that $x$ and $y$ are any elements of some lower section of $P$; then $1_x = 1_y$. If $x < y$, then \eqref{!1} implies that $x \to y = 1_y = 1_x$. If $x = y$, then $x \to y = x * x = 1_x$ by \eqref{1p}. If $y < x$, then (using \eqref{xtoy=}) 
\[
x \to y = \max\{u \ge y\colon [y,u] \cap [y,x] \subseteq \{y\}\} = y,
\] 
as both $u$ and $x$ belong to the chain $[y)$ and are comparable; so, $x \to y < x \le 1_x$. If $x$ and $y$, on the contrary, do not belong to any common lower section, then  $y \nleq x$ and, by \eqref{!1}, $x \to y = 1_y \neq 1_x$. Put together, the four considered cases show that the chosen elements $x$ and $y$ satisfy \eqref{le/to}; so  $(P,\to)$ is left implicative. If, in actuality, $P$ does not fulfill the condition of the corollary, then there is a lower section containing incomparable elements $x$ and $y$. Thus $y \nleq x$ and, therefore $x \to y = 1_y = 1_x$ by \eqref{!1} and the choice of $x$ and $y$. But $x \nleq y$ too; so,  \eqref{le/to} fails and
 $(P,\to)$ is not left implicative. 

(b) Proved similarly: with the exception of the now unnecessary fourth case in the first half of the above proof, it also goes, up to inessential details, for right implicativity. 
\end{proof}

Another result proved with the help of Theorem \ref{nat1y} shows that the natural extension of an SP-poset can be realized, besides \eqref{xtoy=}, also in quite different way.

\begin{theorem}    \label{natisnat}
An operation $\to$ on an SP-poset $(P,*)$ is a natural extension of $*$  if and only if, for all $x,y$,
\begin{equation}   \label{xtoy=2}
x \to y  = \min\{z*y\colon z \in [y,x] \cup \{y\}\} = \min\{z*y\colon z \in ((x] \cup (y]) \cap [y)\}.
\end{equation}
\end{theorem}
\begin{proof}
Let $(P,*)$ be an SP-poset.
Suppose that the first equalitiy in \eqref{xtoy=2} is  fulfilled. Then $\to$ is an extension of $*$: if $y \le x$, then $[y,x] \neq \varnothing$ and $x \to y = \min\{z*y\colon z \in [y,x]\} = x*y$ in view of Lemma \ref{sp-prop}(c). Further, if $y \nleq x$, then $[y,x] = \varnothing$ and $x \to y = \min\{z*y\colon z \in \{y\}\} = 1_y$; see \eqref{1p}. By Theorem \ref{nat1y}, $\to$ is a natural extension of $*$.

Now suppose that $\to$ is a natural extension of $*$. If $y \le x$, then $x \to y = x*y = \min\{z*y\colon z \in [y,x]\}$ (see Lemma \ref{sp-prop}(c)) $= \min\{z*x\colon z \in [y,x] \cup \{y\}\}$. If, on the contrary, $y \nleq x$, then $x \to y = 1_y$ (by Theorem \ref{nat1y}) $= y * y$ (by \eqref{1p}) $= \min\{z*y\colon z \in [y,y]\} = \min\{z*y\colon z \in [y,x] \cup \{y\}\}$. Thus, \eqref{xtoy=2} holds.
\end{proof}

\section{Sectional j-pseudocomplementation revisited} \label{joinex}

\subsection{Normal extensions.}
Let $(P,*)$ be a $*$-poset. In the case when $P$ is an upper semilattice, the operation $*$ can be extended to a total operation $\to$ by the rule
\begin{equation}    \label{vee-ext}
x \to y := (x \vee y) * y.
\end{equation}
This method of extension seems to be first used in \cite{A1} to extend sectional Boolean complementations  and later in \cite{A2} to extend sectional orthocomplementations. In connection with sp-complementation it was mentioned already in \cite{A1}, but was really applied (in lattices) only in \cite{Ch-ext}. (By the way, in lattices \eqref{CLP} also reduces to \eqref{vee-ext}.) We  call an operation $\to$ defined in this way a \emph{join-extension}, in short, a \emph{j-extension of $*$}.

According to Theorem 1 of \cite{C-impl}, an operation $\to$ on an upper semilattice with $1$ is a j-extension of an sp-complementation  if and only if it has the properties (j1)--(j3) mentioned above in Section \ref{related}.2. As we already know, not all arrow posets with $1$ having these properties are ESP-posets. In this section, we shall present another triple of properties possessed only by ESP-posets but equivalent to (j1)--(j3) in upper ESP-semilattices with 1. To this end, we first adjust the extension rule \eqref {vee-ext} to $*$-posets that not necessarily are upper semilattices. 

There are several reasonable ways how to do this. The one which we choose below  was used in \cite[Lemma 3]{C-impl} to recover implication in weak BCK*-algebras from its restriction to upper sections. In the order dual form (for introducting subtraction-like operations) it was applied in \cite{C-subtr}. 
Extensions obtained in this way were called \emph{normal} in \cite{C-subtr}; we shall adopt here this attribute.

Suppose that the sectional operation $*$ on $P$ is antitonic in the first argument and that an operation $\to$ on $P$ is defined by
\begin{equation}    \label{jext}
x \to y := \max\{z*y\colon x,y \le z\}.
\end{equation}
Then $\to$ is an extension of $*$: for $y \le x$, $x \to y = \max\{z*y\colon x \le z\} = x*y$. Moreover, if $P$ is an upper semilattice, then \eqref{jext} reduces to \eqref{vee-ext}. Owing to Lemma \ref{sp-prop}(c), the rule \eqref{jext} can be applied to SP-posets. Observe  that if an SP-poset is finite, then every element $x \to y$ in \eqref{jext} is completely determined by the set of all minimal upper bounds of $x$ and $y$.

It follows from \eqref{jext} that $x \to y$ is defined only if the pair $x,y$ is bounded from above. Thus, a $*$-poset must be  up-directed to have an extension by \eqref{jext}. For an SP-poset this means that it  necessarily has the greatest element (Corollary \ref{p*p}). Generally, the directedness condition is not sufficient for existence of all the necessary maxima.

\begin{example}    \label{notjext}
The subset $\{z*b\colon a,b \le z\}$ of the bounded SP-poset of Example \ref{sp-e} does not have a greatest element and therefore $a \to b$ does not exist.
\deq
\end{example}

\begin{definition}    \label{norm}
The operation $\to$ on a $*$-poset $(P,*)$ obtained as in \eqref{jext} is called a \emph{normal extension of $*$}. The arrow poset $(P,\to)$ is then a  \emph{normal extension of $(P,*)$}. An esp-complementation (an ESP-poset) is said to be \emph{normal} if it is a normal extension of an sp-complementation (resp., of an SP-poset). 
\deq
\end{definition}

For instance, the operation $\to$ in Example \ref{rpc-e} is a normal extension of the operation $*$ (rp-complementation is antitonic in its first argument).  Example \ref{notjext} actually illustrates the following general result.

\begin{theorem} \label{jext-fin}
A finite SP-poset admits a normal extension if and only it is an upper semilattice.
\end{theorem}
\begin{proof}
Let $(P,*)$ be a finite SP-poset. As noticed above, if $P$ is an upper semilattice, then the j-extension of $*$ is normal. For the converse, assume on the contrary that $P$ is not an upper semilattice; of interest is only the non-trivial case when $P$ is still up-directed. 
 
Let $a$ and $b$ are elements of $P$ having not a join (hence, they are incomparable). Then they have at least two minimal common upper bounds, say, $c$ and $d$. Therefore, $a$ and $b$ are lower bounds of $c$ and $d$; without loss of generality, one may assume that $b$ is even a maximal lower bound. Then $b = c \wedge_b d$ by \eqref{mlb}, whence $c \le d*b$ by \eqref{=x*y=}. Now suppose that, in spite of the above assumption, there nevertheless is  a normal extension  $to$ of the operation $*$\,. Then $ a \to b = e*b \ge c*b, d*b$ for some upper bound $e$ of $a,b$.  Choose any lower bound $f \in [b)$ of $c$ and $e$; thus $f \le c \le d*b \le e*b \le f*b$ (Lemma \ref{sp-prop}(c)). Hence,  $f \le b$  by Lemma \ref{sp-prop}(j) and therefore $f = b$. So, $e \wedge_b c = b$ and, by \eqref{=x*y=},  $e \le c*b$, whence $e \le e*b$ and, further, $e \le b$ by Lemma \ref{sp-prop}(j). But then $b \ge a$, a contradiction. Thus,  the operation $*$ actually has not any normal extension.
\end{proof}

A normal extension of an sp-poset which is not an upper semilattice is presented in the subsequent example.

\begin{example}    \label{jext-e1}
Consider a poset consisting of an infinite descending chain $c_1 > c_2 > c_3 > \ldots$ and two  incomparable elements $a,b$ such that $a,b < c_i$ for all $i$. Thus, $c_1$ is the greatest element. This poset is sc-pseudocomplemented: for $j = 1,2,3, \ldots$ and $k > j$, \par
$a * a = 1$,  $b * b = 1$, $c_j * a = a$, $c_j * b = b$,  $c_j * c_k = c_k$, $c_j * c_j = 1$. \\
It has a normal extension where, for $j = 1,2,3, \ldots$, $i < j$, \par
$a \to b = b$, $a \to c_j = 1$,
$b \to a = a$, $b \to c_j = 1$,
$c_j \to c_i = 1$.
\deq
\end{example}

\begin{theorem} \label{natimplext}
The normal extension of an SP-poset is implicative.
\end{theorem}
\begin{proof}
Suppose that $(P,*)$ is an  SP-poset and that $\to$ is a normal extension of the operation $*$. As noted above, then $P$ has a greatest element. The equivalence \eqref{le/to} therefore takes in $(P,\to,1)$ the form
\begin{equation}    \label{le/to1}
x \le y \tiff x \to y = 1.
\end{equation}
It is fulfilled: if $x \le y$, then $x,y \le y$ and, further, $1 = y * y \le x \to y \le 1$ by \eqref{1p}, respectively, (\ref{jext}), and if, conversely, $x \to y = 1$, then $1 = z*y$ for some $z$ with  $x,y \le z$ by \eqref{jext}, whence $z = y$ (see Lemma \ref{sp-prop}(l) and \eqref{1p}) and, further,  $x \le y$.
\end{proof}

\subsection{Axiomatization of normal ESP-posets.} Our next aim is to show that the class \ESP$_\nrm$ of  normal ESP-posets is
characterized by four axioms
\begin{enumerate}
\thitem{nrm0}
$y \le x \to y$,
\thitem{nrm1}
if $x \le y \to z$, then $y \le x \to z$,
\thitem{nrm2}
if $y \le x \le x \to y$, then $x \le y$,
\thitem{nrm3}
if $z$ is a maximal lower bound of $x$ and $y$, then $x \le y \to z$,
\end{enumerate}
the axiom (nrm0) being superfluous if existence of the greatest element of a poset is explicitly assumed. Thus, an esp-complementation is normal if and only if it satisfies the axioms (nrm0) and (nrm1). By the way, then (nrm0) allows one to consider the element $x \to y$ of $P$ merely as a \emph{normal pseudocomplement of $x$ in $[y)$}.

\begin{lemma} \label{jext-prop}
An arrow poset $(P,\to)$ satisfying \textup{(nrm1)} has the properties
\begin{myenum}
\item
$x \le (x \to y) \to y$
\item
if $x \le y$, then $y \to z \le x \to z$,
\item
$(x \to y) \to y) \to y = x \to y$.
\end{myenum}
If it satisfies also \textup{(nrm0)}, then
\begin{enumerate}
\thitem{d}
$x \le y \to y$.
\end{enumerate}
If $(P,\to)$  has a greatest element $1$  and satisfies \textup{(nrm1), (nrm3)}, then
\begin{enumerate}
\thitem{e}
$x \to x = 1$,
\thitem{f}
$y \le x \to y$,
\thitem{g}
$y \le (x \to y) \to y$,
\thitem{h}
$x \to 1 = 1$,
\end{enumerate}
and if also \textup{(nrm2)} is fulfilled, then, moreover,
\begin{enumerate}
\thitem{i}
$1 \to x = x$,
\thitem{j}
$x \le y$ if and only if $x \to y = 1$.
\end{enumerate}
\end{lemma}
\begin{proof}
(a) Trivially by (nrm1). \par
(b) If $x \le y$, then $x \le (y \to z) \to z$ by (a); then (nrm1) applies. \par
(c) By (a) and (b), $((x \to y) \to y) \to y \le x \to y$. The reverse inequality is a particular case of (a). \par
(d) Evidently. \par
(e) By (nrm3), $1 \le x \to x \le 1$. \par
(f) By (nrm1), as $x \le 1 = y \to y$ by (e). \par 
(g) A particular case of (f).   \par
(h) By (f), $1 \le x \to 1 \le 1$. \par 
(i) From (f), $x \le 1 \to x$. By (a), $1 \to x \le 1 \le (1 \to x) \to x$. Then $1 \to x \le x$ by (nrm2). \par
(j) Using (i) and (nrm1), $x \le y$ iff $x \le 1 \to y$ iff $1 \le x \to y \le 1$.
\end{proof}

\begin{corollary}    \label{123}
In any arrow poset $(P,\to)$, \textup{(nrm1)} is equivalent to the conjunction of
items \textup{(a)} and \textup{(b)} of the lemma.
\end{corollary}
\begin{proof}
By the lemma, (nrm1) implies (a) and (b). Conversely, if  (a) and (b) hold and $x \le y \to z$, then $y \le (y \to z) \to z \le x \to z$.
\end{proof}

\begin{corollary}    \label{top}
An arrow poset $(P,\to)$ satisfying \textup{(nrm1)} and \textup{(nrm3)} has a greatest element if and only if it satisfies also \textup{(nrm0)}.
\end{corollary}
\begin{proof}
By Lemma \ref{123}(d,e,f).
\end{proof}

 Now we can prove the main result of the subsection.

\begin{theorem} \label{jext-ax}
An arrow poset is a normal ESP-poset if and only if it satisfies the axioms \textup{(nrm0)--(nrm3)}. The axiom \textup{(nrm0)} may be replaced by a requirement that the poset has a greatest element.
\end{theorem}
\begin{proof}
Let $(P,\to )$ be an arrow poset. Suppose that (nrm0)--(nrm3) are fulfilled in it;  then it satisfies (esp2), (esp3) and, by Lemma \ref{jext-prop}(b), also (esp1), and is therefore an extension of an SP-poset $(P,*)$. For given $x,y$, Lemma \ref{jext-prop}(a,g,c) provides an element $z_0 := (x \to y) \to y$ such that $x,y \le z_0$ and $x \to y = z_0 \to y = z_0*y$. On the other hand, if $x,y \le z$ for some $z$, then $z*y = z \to y \le x \to y$ by Lemma \ref{jext-prop}(b). Thus, \eqref{jext} holds and, consequently, the extension $(P,\to)$ is normal.

This proves that the condition in the first statement is sufficient. To prove that it is also necessary, suppose that $(P,\to)$ is a normal extension of an SP-poset $(P,*)$. As it was noted above,  then $P$ has a greatest element. Moreover, for all $x,y,v$,
\begin{equation}    \label{jext!}
v \le x \to y\ \tiff \mb{there are} u \ge v \tand  z \ge x  \mb{such that} u \wedge_y z = y.
\end{equation}
Indeed, suppose that $v \le x \to y$. By \eqref{jext}, this is the case if and only if $v \le z*y$ with some $z \ge x,y$, i.e.\ (by \eqref{x*y=}), if and only if $v \le u$ for some $u$ such that $u \wedge_y z = y$ and $z \ge x,y$. But the restriction $z \ge y $ here may be omitted as it is implicit in the preceding equality.

Now we can derive the four axioms.

(nrm0)
Put $v = y$ in the equivalence \eqref{jext!}. Then its right-hand side holds for all $x$ and $y$: with $u = y$ and $z = 1$, the condition $u \wedge_y z = y$ is fulfilled. So, its left-hand side also holds.

(nrm1)
By \eqref{jext!}, $x \le y \to z$ iff $x' \wedge_z y' = z$ for some $x' \ge x$ and $y' \ge y$. Likewise $y \le x \to z$ if and only if $y' \wedge_z x' = z$ for some $y' \ge y$ and $x' \ge x$. So, the inequalities in (nrm1) are even equivalent.

(nrm2)
By \eqref{jext!}, if $x \le x \to y$, then $x' \wedge_y x'' = y$ for some $x' \ge x$ and $x'' \ge x$, whence $x'\perp_y x''$. If also $y \le x$, it follows that $x = y$.

(nrm3)
Suppose that $z$ is a maximal lower bound of $x$ and $y$. Then $x \wedge_z y = z$ (see \eqref{mlb}), and $x \le y * z$ immediately by the definition of sp-complementation \eqref{x*y=}. So $x \le y \to z$.

The other statement of the theorem  is immediate from Corollary \ref{top}.
\end{proof}

One more characterization of normal esp-complementation is immediate from the proof of the theorem.

\begin{proposition}	\label{jext-ch}
An operation $\to$ on a poset is a normal esp-complementation if and only if it satisfies \eqref{jext!}.
\end{proposition}

We have also came to the desired conclusion (see the beginning of the section). First, the following are quivalent in an upper semilattice with a greatest element: 
\par -- an operation $\to$ satisfies the axioms (j1)--(j3),  
\par -- it is a join-extended sp-complementation,   
\par -- it satisfies the axioms (nrm1)--(nrm3). 
\\ 
Further, we now know that any arrow poset with 1 satisfying (nrm1)--(nrm3) is an ESP-poset. Therefore, just the operations satisfying these axioms are those which deserve to be called sectional j-pseudocomplementations. However, we are not going to alter here  the present meaning of the term.

The two axiom systems are equivalent also in lower arrow semilattices. Notice that the  proof below goes also for upper semilattices..

\begin{proposition}    \label{j'=j}
An operation $\to$ on a lower semilattice with unit is a normal esp-complementation if and only if it is a sectional j-pseudocomplemen\-tation.
\end{proposition}
\begin{proof}
The axioms (j1) and (nrm1) coincide, while (j3) and (nrm3) both amount to (esp$^\wedge$2) in lower semilattices. Next, (nrm2) is contained in (j2), but (j2) is fullfilled in every normal ESP-poset that is a lower semilattice: by \eqref{jext!}, if $x \le x \to y$, then $x' \wedge_y x'' = y$ for some $x' \ge x$ and $x'' \ge x$. So $y = x'\wedge x''$ in view of Proposition \ref{glb}, whence $x \le y$.
\end{proof}

Like general implicative ESP-algebras  (Remark \ref{notQV} in Section \ref{ext}), also the class of all normal implicative ESP-algebras (see Theorem \ref{natimplext}) is not closed under forming subalgebras, as the following example shows, and therefore is not a quasivariety. 

\begin{example}
Consider the poset $P$ containing elements $0,a.b,c,1$ and formed by two maximal chains $0 < a < b < 1$ and $0 < a < c < 1$. Equipped with an operation $\to$ having the table
\[
\begin{array}{c|ccccc}
\to &0&a&b&c&1 \\
\hline
0   &1&1&1&1&1 \\
a   &0&1&1&1&1 \\
b   &0&c&1&c&1 \\
c   &0&b&b&1&1 \\
1   &0&a&b&c&1
\end{array}
\]
it yields an implicative ESP-algebra that satisfies the axioms (nrm1)--(nrm3). Its subalgebra $Q$ consisting of elements $0,b,c,1$ is not an ESP-poset at all: the sp-complement of $b$ in $[0)$ in $Q$ equals $c$.
\deq
\end{example}

\subsection{Normal ESP-semilattices and lattices.}
In contrast, normal ESP-meet semilattices do form a subvariety of \ESP$^\wedge$. By the help of Remark \ref{sp=wr} and Corollary \ref{j'=j}, this could be inferred from Proposition 23 of \cite{C-impl}; however, the proof of that proposition is not accurate. The theorem below is even a strengthening of the proposition (recall the note preceding Theorem \ref{ESP^}).

\begin{theorem}  \label{nESP^}
The class \ESP$^\wedge_\nrm$ of meet semilattices with normal esp-comple\-mentation is a subvariety of \ESP$^\wedge$ determined by the identities
 \\[2pt]
\textup{(nrm$^\wedge$0)}
$y \le x \to y$, \\
\textup{(nrm$^\wedge$1)}
$x \le (x \to y) \to y$, \\
\textup{(nrm$^\wedge$2)}
$x \to z \le (x \wedge y) \to z$. 
\end{theorem}
\begin{proof}
Let $(P,\wedge,\to)$ be a meet ESP-semilattice. The axiom (nrm0) coincides with (nrm$^\wedge$0). By Corollary \ref{123}, the axiom (nrm1) is equivalent in $(P,\to)$ to the conjunction of (nrm$^\wedge$1) and (nrm$^\wedge$2) as (nrm$^\wedge$2) is another form of (nrm1). In addition, (nrm2) and (nrm3) coincide with (esp2) and (esp3), respectively. Thus, $\to$ is a normal esp-complementation if and only if (nrm$^\wedge$0)--(nrm$^\wedge$2) are fulfilled in $(P,\wedge,\to)$. 
\end{proof}

In the context of normal extensions, the \ESP$^\wedge$-axiom (esp$^\wedge$1) can be simplified.

\begin{corollary}
The variety \ESP$^\wedge_\nrm$ is characterized by semilattice axioms, the axioms \textup{(nrm$^\wedge$0)}--\textup{(nrm$^\wedge$2)} and  
\\[2pt]
\textup{(nrm$^\wedge$3)} \ $x \wedge (x \to y) \le y$,  \\ 
\textup{(nrm$^\wedge$4)} \ $x \le y \to (x \wedge y)$. 
\end{corollary}
\begin{proof}
As (nrm$^\wedge$4) coincides with  (esp$^\wedge$2), it remains to show that the axiom (esp$^\wedge$1)  of ESP$^\wedge_{nrm}$ (inherited from ESP$^\wedge$) can be replaced by the simpler axiom (nrm$^\wedge$3). Evidently,  (nrm$^\wedge$3) implies (esp$^\wedge$1). Conversely, (nrm$^\wedge$3) holds in ESP$^\wedge_{nrm}$: put $v := x \wedge (x \to y)$; then
 $v \le x$ and $v \le x \to y$, whence $x \to y \le v \to y$ by (nrm$^\wedge$2) and, further, $v \le y$ by (j2) (which is fulfilled in \ESP$^\wedge_\nrm$ by Proposition \ref{j'=j}). 
\end{proof}

In fact, (j2) itself  is an easy consequence of (nrm$^\wedge$3). In its turn, the identity    (nrm$^\wedge$3) can be shorthened, in view of (nrm$^\wedge$0), to 
\\[3pt]
\textup{(nrm$^\wedge$3')} \ $x \wedge (x \to y) = x \wedge y$.
\\[3pt]
Moreover, in lower semilattices also the equivalence \eqref{jext!} can be refined:
\begin{equation}    \label{jext!wedge}
v \le x \to y\ \tiff \mb{there are} u \ge v \tand  z \ge x  \mb{such that} u \wedge z = y;
\end{equation}
see Lemma \eqref{glb}. 

The equivalence \eqref{jext!} can be simplified also in upper semilattices. But recall that in upper semilattices (and  lattices) the extension rule \eqref{jext} reduces to \eqref{vee-ext}; so referencing to normality becomes contingent in a sense, or outer, in these cases. To unbind the notation and terminology from this contingency, we return to the terms  `join-extension' and `j-extension' introduced at the beginning of the section, and write \ESP$^\vee_\j$, \ESP$^{\wedge\vee}_\j$ instead of \ESP$^\vee_\nrm$, \ESP$^{\wedge\vee}_\nrm$ for the class of join semilattices, resp.,  lattices with join-extended (thus, normally extended) sp-complementation.  

\begin{proposition}[\protect{\cite[Lemma 1]{C-impl}}]  \label{v-Jext!}
In an upper semilattice, an operation $\to$ is a j-extended sp-complementation if and only if, for all $x,y$  and $v$,
\begin{equation} 	\label{v-jext!}
v \le x \to y \ \tiff \ y \mb{is  a meet of} v \vee y \tand x \vee y.
\end{equation}
\end{proposition}
\begin{proof}
Let $(P,*)$ be an upper ESP-semilattice, and let $\to$ be a extension of $*$. Then
\begin{multline*}  
v \le (x \vee y) * y  \tiff  v \vee y \le (x \vee y) * y  \\
\tiff  (v \vee y) \wedge _y (x \vee y) = y \tiff (v \vee y) \wedge (x \vee y) = y.
\end{multline*}
by Lemma \ref{sp-prop}(a),  \eqref{in*} and Lemma \ref{glb}. Together with  \eqref{vee-ext}, this leads to \eqref{v-jext!}, and together with  \eqref{v-jext!}, to \eqref{vee-ext}.
\end{proof}

Dualizing the traditional notion of nearlattice (see, e.g., \cite{C-subtr} and references therein), we assume the following definition. 

\begin{definition}
An upper semilattice where every pair of elements  bounded from below has a greatest lower bound is called a \emph{nearlattice}.
\deq
\end{definition}

In view of  Proposition \ref{glb}, this definition is equivalent to that used in \cite{ChKo}. Evidently, every finite upper semilattice is a nearlattice. The finiteness condition here is essential.

Finally, an upper semilattice is a nearlattice if and only if the meet $(x \vee y) \wedge (y \vee z)$ exists in it for all $x,y,z$. Thereby, the meet in the right-hand part of \eqref{v-jext!} is, in a nearlattice, always defined. Owing to this, j-extended sp-complementations on nearlattices  can be characterized by existentially interpreted partial identities  (``both sides of the equation are defined, and they are equal'') .
An inequality $a \le b$ in the following theorem is considered as a shorthand for the equality $a \vee b = b$.  

\begin{theorem}
Let $(P,\to)$ be an upper arrow semilattice. It is an j-extended SP-nearlattice if and only if 
the following conditions are satisfied:
\begin{enumerate}
\thitem{j$^{\wedge\vee}$1} $(x \to y) \wedge (x \vee y) = y$,
\thitem{j$^{\wedge\vee}$2} $z \le x \to ((z \vee y) \wedge (x \vee y))$,
\end{enumerate}
where the involved meets are assumed to be defined.
\end{theorem}
\begin{proof}
	Assume that $P$ is a nearlattice and  $\to$ is a j-extended sp-complementation on $P$. The equivalence \eqref{v-jext!} yields  identities
\[
y \le x \to y \ \tand \   ((x \to y) \vee y) \wedge (x \vee y) = y,
\]
which in conjunction imply (j$^{\wedge\vee}$1). The axiom  (j$^{\wedge\vee}$2) also follows from the \textsl{if} part of the same equivalence by the substitution of $z$ for $v$ and  the meet $(z \vee y) \wedge (x \vee y)$ (which is always defined in a nearlattice) for $y$: it is easily seen that 
\[
z \vee ((z \vee y) \wedge (x \vee y)) = z \vee y\ \mbox{ and likewise }\ x \vee ((z \vee y) \wedge (x \vee y)) = x \vee y.
\]
Thus a j-extended sp-complementation on $P$ satisfies the two axioms.

Now assume that, conversely, the operation $\to$ on a semilattice $P$ satisfies the axioms. The \textsl{if} part of \eqref{v-jext!} follows from (j$^{\wedge\vee}$2).
 To prove its \textsl{only if} part, suppose that $z \le x \to y$ and recall that the axioms are existential identities, i.e., the meets $(z \vee y) \wedge (x \vee y)$ and $(x \to y) \wedge (x \vee y)$ are defined in $P$ for all $x,y,z$. From (j$^{\wedge\vee}$1), $y \le x \to y$; so  $z \vee y \le x \to y$ and, further, $y \le (z \vee y) \wedge (x \vee y)\le (x \to y) \wedge (x \vee y) = y$ by (j$^{\wedge\vee}$1). Therefore, the supposition implies that $(z \vee y) \wedge (x \vee y) = y$. Thus \eqref{v-jext!} holds and $\to$ is a j-extended sp-complementation. In addition, $P$ is a nearlattice: if $c \le a,b$ in $P$, then $a = a \vee c$, $b = b \vee c$, and the meet $a \wedge b$ exists.
\end{proof}

These axioms do not work under the conditional interpretations of equations ``if both sides of the equation are defined, then they are equal''  and ``if one side of the equation is defined, then the other also is and they are equal''.

Of course, every lattice is a nearlattice; so the following conclusion (where $\wedge$ is now a lattice meet) from the theorem is evident. 

\begin{corollary}	\label{char-jwedvee}
The class  \ESP$^{\wedge\vee}_\j$ is a variety characterized by the lattice axioms and the identities \textup{(j$^{\wedge\vee}$1), (j$^{\wedge\vee}$2)}.
\end{corollary}

\begin{remark}
Nearlattices with j-extended sp-complementation were first axiomatized in \cite[Theorem 5]{ChKo}; the so called almost subtractive nearlattices, which were introduced and studied in \cite{C-subtr},   are essentially the order duals of these (to take an example,  Proposition \ref{v-Jext!} above corresponds to Corollary 3.5 in \cite{C-subtr}). Nearsemilattices with pseudocomplemented sections and their j-extensions have recently been considered also in  \cite{ChL2}.  The result that the  lattices with j-extended sp-complementation form a variety was first  stated in \cite[Theorem 2]{Ch-ext}.  The pair of axioms (j$^{\wedge\vee}$1), (j$^{\wedge\vee}$2) for them is evidently equivalent to the apparently stronger one presented in Theorem 2.9 of \cite{ChLP} (these lattices are called there simply sectionally pseudocomplemented). The axiom (j$^{\wedge\vee}$2) will be further simplified in Section \ref{another}.3.
\deq
\end{remark}

\section{Generalized natural extensions}   \label{another}

\subsection{Local subset selection functions}

The definition \eqref{xtoy=} can be given another form. We first rewrite the definition \eqref{x*y=} as follows, where, in particular, the partial operation $\wedge_y$ is eliminated:
\begin{equation}    \label{.x*y=}
x*y = \max\{u\colon (u] \cap (x] \cap [y) = \{y\}\}.
\end{equation}
The operation $\to$ in \eqref{xtoy=} admits a structurally similar characteristic:
\begin{equation}  \label{.xtoy=}
x \to y =  \max\{u\colon (u] \cap ((x] \cup (y]) \cap [y) = \{y\}\} .
\end{equation}
Indeed,
\begin{align*}   
\max\{u \ge y\colon u \perp_y x\} &= \max\{u \ge u\colon (u] \cap (x] \cap [y) \subseteq \{y\}\} \\ 
 &= \max\{u \ge y\colon (u] \cap ((x] \cup (y]) \cap [y) \subseteq \{y\}\} \\
 &= \max\{u\colon (u] \cap ((x] \cup (y]) \cap [y) = \{y\}\} . 
\end{align*}
Formally, the equality \eqref{.xtoy=} is obtained from  \eqref{.x*y=} by replacing the component $(x]$ in the defining part of the former by $(x] \cup (y]$, but more general reasonable replacements are also possible. 

Let $P$ be any poset. We associate with every pair of elements $x,y$ of $P$ a down set $I(x,y)$ so that

\begin{enumerate}
\item[(I0)] $x,y \in I(x,y)$
\item[(I1)]  $I(x,y) = I(y,x)$,
\item[(I2)] if $y \le x$, then $I(x,y) = (x]$,
\item[(I3)] if $x \le x'$, then $I(x,y) \subseteq I(x',y)$;
\end{enumerate}
the order ideals $(x] \cup (y]$ appearing in \eqref{.xtoy=} is an example of such a family of sets $I(x,y)$. Notice that $I(x,x) = (x]$ and that  in an upper semilatice $I(x,y) \subseteq I(x \vee y, y) = (x \vee y]$. More generally, it follows from (I3) and (I2) that
\begin{enumerate}
\item[(I4)] if $x,y \le z$, then $I(x,y) \subseteq (z]$.
\end{enumerate}
So, always
\begin{enumerate}
\item[(I5)] $(x] \cup (y] \subseteq I(x,y) \subseteq \bigcap((z]\colon x,y \le z)$.
\end{enumerate}
A function $I$ with the described properties thus extends the embedding $x \mapsto (x]$ of $P$ into the set of down sets of $P$. We consider every subset $I(x,y)$ as a selected substitute for the (possibly unexistent) lower section $(x \vee y]$, and call any function $I$ having the properties (I0)--(I3) a \emph{local subset selection}.  Besides  the order ideals, on which  \eqref{.xtoy=} rests, good candidates for the sets $I(x,y)$ are poset ideals of any other kind generated by $\{x,y\}$. More generally, even any standard completion of a poset $P$ (i.e., a  closure system of its down sets containing all principal ideals; see, e.g., \cite[Section 1]{E} for examples) gives rise, in the same way, to a local subset selection. We shall touch the local subset selection induced by Frink ideals in the final subsection.

\subsection{Extension with the help of local selections}

In this subsection, let $I$ stand for some local subset selection on an arrow poset $(P,\to)$. 

Consider a generalized form of \eqref{.xtoy=}:
\begin{align}     \label{Ixtoy=}
x \to y  &:= \max\{u\colon (u] \cap I(x,y) \cap [y) = \{y\}\} \\
& ~= \max\{u \ge y\colon (u] \cap I(x,y) \cap [y) \subseteq \{y\}\}.	\notag
\end{align}
The operation $\to$ here is still an  esp-complementation: when $y \le x$, (I2) implies that $x \to y$ is a pseudocomplement of $x$ in $[y)$; cf.\ \eqref{.x*y=}. 

\begin{definition}
The operation $\to$ is called an \emph{$I$-natural} esp-complemen\-tation if it satisfies \eqref{Ixtoy=}. The arrow poset itself is then called an $I$-natural ESP-poset.
\deq
\end{definition}

Observe that, for any local selections $I'$ and  all $x,y$,
\begin{align*}
I(x,y) \subseteq I'(x,y) & \Rightarrow (u] \cap I(x,y) \cap [y) \subseteq (u] \cap I'(x,y) \cap [y) \\
& \Rightarrow \tif (u] \cap I'(x,y) \cap [y) \subseteq \{y\}, \tthen  (u] \cap I(x,y) \cap [y) \subseteq \{y\} \\
& \Rightarrow \max\{u \ge y\colon (u] \cap I'(x,y) \cap [y) \subseteq \{y\}\} \\
     &  {\qquad\qquad}\hfill \le \max\{u \ge y\colon (u] \cap I(x,y) \cap [y) \subseteq \{y\}\} \\
& \Rightarrow x \to 'y \le x \to y
\end{align*}
(provided that both involved maxima exist), where $\to'$ is the esp-complementation induced by $I'$. 

The definition \eqref{Ixtoy=} of the operation $\to$ can be transformed into a more convenient  implicit definition consisting of three conditions
\begin{multline}    \label{I=xtoy=}
\textup{(i) } y \le x \to y, \ \textup{(ii) } (x \to y] \cap I(x,y) \cap [y) \subseteq \{y\}, \\
\textup{(iii) } \tif y \le u\tand (u] \cap I(x,y) \cap [y) \subseteq \{y\}, \tthen u \le x \to y
\end{multline}
 on $\to$. Just as in \eqref{=xtoy=}, the first of them is a consequence of the last one. The properties of $\to$ listed in the following lemma either strengthen or supplement those given in Proposition \ref{esp-prop}.

\begin{lemma}    \label{Inat-prop}
Independently of the choice of the function $I$, an $I$-natural esp-complementation has the following properties:
\begin{myenum}
\item 
$y \le x \to y$,
\item 
$y \le (x \to y) \to y$,
\item 
if $x \le y$, then $x \to y = 1_x = 1_y$,
\item 
if $x \le y$, then $y \to z \le x \to z$,
\item 
if $z \le y$ and $x \le y \to z$, then $y \le x \to z$,
\item  
$x \to y \le 1_y$,
\item
$x \to 1_y =1_y$.
\end{myenum}
\end{lemma}
\begin{proof}
(a) From \eqref{I=xtoy=}(i).

(b) From (a).

(c) If $x \le y$, then $I(x,y) = (y]$ (by (I2)) and $x \to y = \max\{u\colon (u] \cap \{y\} = \{y\}\} = 1_y = 1_x$, see Lemma \ref{esp-prop}(i).

(d) Suppose that $x \le y$. By (I3) and \eqref{I=xtoy=}(ii),
\[
(y \to z] \cap I(x,z) \cap [z) \subseteq (y \to z] \cap I(y,z) \cap [z) \subseteq \{z\}
 \]
and, further, $y \to z \le x \to z$ by \eqref{I=xtoy=}(iii,i).

(e) Suppose that $z \le y$. If $x \le y \to z$, then $(y \to z) \to z \le x \to z$ by (d), whence $y \le x \to z$ by Lemma \ref{esp-prop}(b).

(f)  By (a) and  Lemma \ref{esp-prop}(g).

(g) By (a), $1_y \le x \to 1_y$. The reverse follows from (f), as $\max [1_y) = 1_y$.
\end{proof}

Under an certain assumption on the operation $\to$,  the condition (ii) in \eqref{I=xtoy=} can be replaced by a set of simpler ones: the function $I$ can be excluded. However, some preliminary work for this has to be done.

\begin{lemma}	\label{simplI}
For all $u,x,y$,
\begin{myenum}
\item
\(  (u] \cap I(x,y) \cap [y) \subseteq \{y\} \ \tiff \ u \perp_y z  \tforall z \in I(x,y) \), 
\item
\(  (u] \cap I(x,y) \cap [y) = \{y\} \ \tiff \ u \wedge_y z  = y \tforall z \in I(x,y) \cap [y) \).
\end{myenum}
\end{lemma}
\begin{proof}
(a)
Since $I(x,y) = \bigcup((z]\colon z \in I(x,y))$, 
\begin{align*}
(u] \cap I(x,y) \cap [y) \subseteq \{y\}
& \tiff [y,u] \cap \bigcup([y,z]\colon z \in I(x,y)) \subseteq \{y\} \\
& \tiff \bigcup([y,u] \cap [y,z])\colon z \in I(x,y)) \subseteq \{y\} \\
& \tiff \tforall  z \in I(x,y), [y,u] \cap [y,z] \subseteq \{y\} \\ 
& \tiff \tforall  z \in I(x,y),  u\perp_y z. 
\end{align*}

(b) Since $u \perp_y z$ whenever $y \nleq z$, it follows that the right part of (a) holds if and only if  $u \perp_y z$ for all $z \in I(x,y) \cap [y)$, i.e., if and only if ($u \perp_y z$ and $y \le z$) for all $z \in I(x,y) \cap [y)$.Further
, the left part of (b) is equivalent to the conjunction of an inequality $y \le u$ and the left part of (a), and the same holds for both right parts. So (b) follows from (a). 
\end{proof}

The idea of the following definition is borrowed from \cite[Definition 2.1]{ChLP}.

\begin{definition}
An esp-complementation $\to$ is \emph{strong} if it satisfies the condition
\\[2pt]
(S) $x \le (x \to y) \to y$
\\[2pt]
from Lemma \ref{jext-prop}(a). A \emph{strong} ESP-poset is one with a strong esp-complemen\-tation. ~
\deq
\end{definition}

For instance, the ESP-poset in Example \ref{sp-e1} is strong. Also, the normal esp-complementation of any SP-poset is strong (Lemma \ref{jext-prop}(a)). Example \ref{neither} below shows that, for a certain $I$, there are posets, where an $I$-natural esp-complementation is not strong.

\begin{proposition}	\label{str-nrm}
Every strong $I$-natural ESP-poset is normal.
\end{proposition}
\begin{proof}
Let $(P,\to)$ be an ESP-poset satisfying both suppositions.
By Lemma \ref{Inat-prop}(a,d) and Corollary \ref{123}, it then satisfies (nrm0) and (nrm1). In view of (esp2) and (esp3), it satisfies also (nrm2) and (nrm3). By Theorem \ref{jext-ax}, $(P,\to)$ is  normal.
\end{proof}

Therefore a strong $I$-natural ESP-poset has a greatest element and is implicative, and if it is finite, then it is an upper semilattice. Example \ref{neither} shows that strongness condition here is essential.

\begin{theorem}	\label{sInat}
An up-directed arrow poset is an I-natural ESP-poset if and only if it satisfies the axioms \textup{(nat1), (nat2)} and
\\[2mm]
\textup{(nat$_I$3)}
if $z \le x$ and $x \perp_z  w  \tforall w \in I(y,z )$, then $x \le y \to z$.
\end{theorem}
\begin{proof}
Let $(P,\to)$ be an  arrow poset in which (S) is fulfilled. If the operation $\to$ is an I-natural esp-complementation, then it satisfies \eqref{I=xtoy=}(iii), i.e., the axiom (nat$_I$3) (see Lemma \ref{simplI}(a)), and, by Lemma \ref{Inat-prop}, also (nat1) and (nat2).

Now suppose that (nat1), (nat2) and (nat$_I$3) are fulfilled in $(P,\to)$. Then (nat$_I$3) gives us \eqref{I=xtoy=}(i). To prove \eqref{I=xtoy=}(ii), we, in view of Lemma \ref{simplI}(a), have to show that $(x \to y) \perp_y w  \tforall w \in I(x,y)$. So, let $w \in I(x,y)$. 
By (I5), then $w \le z$ for all $z \ge x,y$. Hence, $x \to y \le z \to y$ by (nat1). It further follows for every $v$ with $y \le v \le x \to y,w$ that $v \le z \to y$ and $z \to y \le v \to y$ (again by (nat1), as $v \le w \le z$) and $v \le y$ (by (nat2)), as needed. So, $(x \to y) \perp_y w$, as required. Finally, \eqref{I=xtoy=}(iii) follows from (nat$_I$3), and \eqref{I=xtoy=} holds.
\end{proof}

The axiom (nat$_I$3) reduces to (nat3) when $I(y,z) = (y] \cup (z]$. The following corollary to the theorem is immediate. 

\begin{corollary}	\label{iso}
Suppose that local selection $I'$ is such that $I(x,y) \subseteq I'(x,y)$ for all $x,y$. 
Then every up-directed $I$-natural ESP-poset is $I'$-natural.
\end{corollary}

\subsection{Extensions via Frink ideals}

A subset $J$ of $P$ is a \emph{Frink ideal} if $L(U(J)) \subseteq J$ whenever $J$ is a finite subset of $J$; see, e.g., \cite{E}. Let $F$ be the local subset selection that associates with every pair of elements $x,y \in P$ the Frink ideal generated by them:
\[
F(x,y) := L([x) \cap [y)) = \bigcap((z]\colon x,y \in (z])
\]
(therefore, $F(x,y) = P$ if $x$ and $y$ have no common upper bound). So,  \eqref{Ixtoy=} takes the form
\begin{equation}    \label{Fxtoy=}
x \to y := \max\{u\colon (u] \cap L(([x) \cap [y)) \cap [y) = \{y\}\},
\end{equation}
which may be considered  as a correct sp-counterpart of \eqref{CLP}; observe that, for every $u$, the equality (*) from Section \ref{related}.3 (i.e., $(u] \cap L([x) \cap [y)) = (y]$) implies that $(u] \cap L([x) \cap [y)) \cap [y)= \{y\}$.  
In particular,
\[
\mbox{$x \to y = y$ if $x$ and $y$ have no common upper bound}.
\]
It follows from Corollary \ref{iso} and (I5) that every strong $I$-natural ESP-poset is $F$-natural.

A poset with F-natural esp-complementation may be neither implicative nor strong and, hence, may be not normal  (cf. Proposition \ref{str-nrm}) even when it has a greatest element. 

\begin{example}    \label{neither}
In the SP-poset $(P,*)$ from Example \ref{sp-e}, the F-natural esp-com\-plementation $\to$ of $*$ has the table
\[
\begin{array}{c|cccccc}
\to &0&a&b&c&d&1 \\
\hline
0 &1&1&1&1&1&1 \\
a &b&1&1&1&1&1 \\
b &a&1&1&1&1&1 \\
c &0&d&d&1&d&1 \\
d &0&c&c&c&1&1 \\
1 &0&a&b&c&d&1
\end{array}
\]
The extension $(P,\to)$ is not implicative, for $a \to b = b \to a = 1$, and not strong, for $(a \to b) \to b = b$. It differs from the pure extension of $(P,*)$ (Example \ref{sp-e1}), where $a \to b = b$ and $b \to a = a$, and from its natural extension in the sense of Section \ref{elim}, where it should hold that $c \to d = d \to c = 1$.  It differs also from the arrow poset of Example \ref{rpc-e}, which, as we now know, is both an rp-poset and a CLP-poset.
\deq
\end{example}

On the other hand, a normal ESP-poset is not necessarily $F$-natural.

\begin{example}
The normal ESP-poset from Example \ref{jext-e1} is not $F$-natural: there, $F(a,b) \cap [b) = \{b\}$ but $a \to b \neq 1$. 
\end{example}

 $F$-natural ESP-lattices are characterized by simple identities. 

\begin{theorem}
An operation $\to$ on a lattice is an $F$-natural esp-comple\-mentation if and only if it satisfies the axioms
\begin{enumerate}
\thitem{j$^{\wedge\vee}$1}
$(x \to y) \wedge (x \vee y) = y$,
\thitem{j$^{\wedge\vee}$2'}
$z \le x \to(z \wedge (x \vee y))$.
\end{enumerate}
\end{theorem}
\begin{proof}
Let $(P,\wedge,\vee)$ be a lattice. Then $F(x,y) = (x \vee y]$ and \eqref{Fxtoy=} reduces to
\begin{equation}   \label{vee_xtoy=}
x \to y = \max\{u\colon (u] \cap  (x \vee y] \cap [y) = \{y\}\}.
\end{equation}
This equality can be presented as a conjunction of two (universally quantified) statements
\begin{enumerate}
\thitem{i} $(x \to y] \cap (x \vee y] \cap [y) = \{y\}$ and 
\thitem{ii} if $(u] \cap (x \vee y] \cap [y) = \{y\}$, then $u \le x \to y$.
\end{enumerate}
But the statement (i) means that $(x \to y) \wedge_y (x \vee y) = y$; by Proposition \ref{glb}, the latter equality is equivalent to the axiom  (j$^{\wedge\vee}$1).  Similarly, the statement (ii) is equivalent to the following one:
\[	
\tif u \wedge (x \vee y) = y, \tthen u \le x \to y.
\]
It is implied by the axiom (j$^{\wedge\vee}$2'). On the other hand, the axiom follows from it by substituting $z$ for $u$ and $z \wedge (x \vee y)$ for $y$ since  $z \wedge (x \vee (z \wedge (x \vee y))) = z \wedge (x \vee y)$ (notice that $x \vee (z \wedge (x \vee y)) \le x \vee (x \vee y) = x \vee y$\,).  Consequently,  \eqref{vee_xtoy=} is equivalent to the conjunction of  (j$^{\wedge\vee}$1) and (j$^{\wedge\vee}$2').
\end{proof}

On the other side, the equality \eqref{vee_xtoy=}, which fits in the general scheme \eqref{Ixtoy=}, is an expanded form of the definition \eqref{vee-ext} obtained by the help of \eqref{.x*y=} and
describes join-extended sp-complementations: it can be rewritten as
\[
x \to y = \max\{u\colon u \wedge_y (x \vee y) = y\}.
\]
We thus have improved Corollary \ref{char-jwedvee} (and also \cite[Theorem 2.9]{ChLP}):

\begin{corollary}	\label{char-jwedvee2}
The class  \ESP$^{\wedge\vee}_\j$ is a variety characterized by the lattice axioms and the identities \textup{(j$^{\wedge\vee}$1), (j$^{\wedge\vee}$2')}.
\end{corollary}  

Of course, also other concrete choices of $I$ may yield interesting  ESP-comple\-men\-tations.

\section{Additional notes}

Thus, there are various reasonable ways how to extend sectional pseudocomplementation on posets.  At least two  motivations for  considering extended SP-posets, both known in the literature, can be mentioned: \par
 (i) ESP-posets, semilattices and lattices, being total algebras, are more suitable for applying methods of universal algebra for study of their structure and the structure of classes of such algebras, \par
 (ii) in algebraic logic, ESP-posets may be viewed as models of many-valued propositional logics with non-Boolean implication.

We do not discuss the item (ii) here. In conection with (i), it is of importance to know in what extent  an extension of an SP-poset reflects the algebraic structure of the latter. For instance, consider a subset $Q$ of a star poset $(P,*)$ as a  subalgebra of it if,  for all $x,y \in Q$, the pseudocomplement $x*y$ belongs to $Q$ whenever it its defined in $P$. Now, is this the case that this happens if and only if $Q$ is a subalgebra of, say, the natural extension of $(P,*)$? Sufficiency of this condition is evident; the proof that it is necessary, is also straightforward if using \eqref{purerule}. Similar questions should be answered for appropriate notions of congruence, homomorphism etc. of SP-posets, and they can be repeated also for normal and $I$-natural extensions of SP-posets. An additional question arises in these cases: which SP-posets are indeed extendible? This is one possible line of further investigations.
 
Another line: besides the normal and $I$-natural extensions, which were chosen partly to cope with the noticed shortcomings of the approaches reviewed in subsections \ref{related}.2 and \ref{related}.3, there are  also various other regular ways for extending SP-posets. As the first one we consider an extension rule dual in a sense to the rule \eqref{jext}.  

The condition for selecting $z$ in the latter can be given a form $(x] \cup (y] \subseteq (z]$.  This suggests a dual condition $[x) \cap [y) \subseteq [z)$ together with looking for the least value of $z*y$. However, $z$ still should satisfy also the inequality $y \le z$. We thus come to the following  definition for a possible extension of $*$:
\begin{equation}    \label{jext2}
x \to y := \min\{z*y\colon  [x) \cap [y) \subseteq [z) \subseteq [y)\}.
\end{equation}
Again, it it is helpful if the operation $*$ is antitonic in the first argument. For $y \le x$, $x \to y = \min\{z*y\colon  y \le z \le x\} = x*y$, i.e., $\to$ is indeed an extension of $*$.  If $P$ is an upper semilattice, then $x \to y = \min\{z*y\colon y \le z \le x \vee y\} = (x \vee y)*y$, i.e.,  \eqref{jext2} reduces to \eqref{vee-ext}. Thus $\to$ turns out to be a j-extended sp-complementation in an upper SP-semilattice; properties of this operation on arbitrary SP-posets are to be investigated. 

The rule \eqref{xtoy=2}, which also yields natural extensions, can be generalized, like \eqref{Ixtoy=}, by the help of local subset selections:
\begin{equation}	\label{Ixtoy=2}
x \to y := \min\{z*y\colon z \in I(x,y) \cap [y)\};
\end{equation}
when $y \le x$, $I(x,y) = (y]$, and $x \to y$ is a pseudocomplement of $x$ in $[y)$. However, it seems that Theorem \ref{natisnat} does not admit a direct generalization for arbitrary function $I$. If indeed so, then the rule \eqref{Ixtoy=2} gives rise to one more family of classes of extended ESP-posets. But it is worth to note that the rule reduces to \eqref{jext2} when $I = F$: $[x) \cap [y) \subseteq [z) \subseteq [y)$ iff $y \le z \le v$ for all $v \ge x,y$ iff $z \in L([x) \cap [y)) \cap [y)$.

Finally, the operation $\to$ on a lower $*$-semilattce obtained by  the rule
\begin{equation}	\label{mext}
x \to y := x * (x \wedge y),
\end{equation}
 evidently is an extension of the operation $*$\,, and, by analogy with j-extensions, may be called an \emph{m-extension} of $*$ (with `m' for `meet'). This rule was first used in \cite{ChH}.  Defined in another way---by
\begin{equation}	\label{sch}
x \to y := \max\{u\colon u \wedge x = x \wedge y\},
\end{equation}
m-extended meet SP-semilattices  were investigated already in \cite{Sch} (the connection with m-extensions appears there as  equality (2.6)) under the name \emph{semi-Brouwerian} semilattices. The element $x \to y$ defined in this way was called  in \cite{Sch} a weak (relative) pseudocomplement of $x$ with respect  to $y$ (in fact, $\to$ is there an extension of week relative pseudocomplementation of Section 4.3).  A question is how could these two rules be adjusted to arbitrary SP-posets. As to \eqref{mext}, \eqref{jext} suggests a similar rule
\[
x \to y := \min\{x*z\colon z \le x,y\},
\]
but it works well only if  the operation $*$ is antitonic in the second argument. One may also first to rewrite the equality $u \wedge x = x \wedge y$ in  \eqref{sch} as an equivalence
\[
\mbox{for all $z$,\, $u \wedge x = z$ if and only if $x \wedge y = z$},
\]
and then apply Proposition \ref{glb} to get a rule
\[
x \to y := \max\{u\colon \mbox{the pairs $u,x$ and $x,y$ have the same maximal lower bounds}\} 
\]
which makes sense in arbitrary posets. The operation $\to$ defined in this way is indeed an esp-complementation: if $y \le x$, then $y$ is the single maximal lower bound of $x$ and $y$, and the latter rule reduces to \eqref{x*y=}; see \eqref{glb}.

\end{document}